\documentclass[11pt,a4]{article}

\usepackage[utf8]{inputenc}
\usepackage[english,main=russian]{babel}

\usepackage{amsfonts, amsmath, amsthm}
\usepackage{amssymb}
\usepackage{url}
\usepackage{cite}

\usepackage{mathtools}
\usepackage{makecell}
\usepackage{tikz}
\usepackage{dsfont}
\usepackage{secdot}

\usepackage{algorithm}
\usepackage{algorithmic}
\usepackage{float}

\sectiondot{section}
\sectiondot{subsection}
\theoremstyle{definition}

\newcommand{\norm}[1]{\left\lVert#1\right\rVert}

\newtheorem{lemma}{Лемма}

\newtheorem{remark}{Замечание}
\newtheorem{statement}{Утверждение}
\makeatletter
\renewcommand{\ALG@name}{Алгоритм}

\allowdisplaybreaks

\begin{document}

\textbf{УДК} 519.85

\begin{center}
\textbf{Ускоренные методы для седловых задач\footnote{Исследование (в § 1, 2) выполнено в рамках Программы фундаментальных исследований НИУ ВШЭ и финансировалось в рамках господдержки ведущих университетов Российской Федерации "5-100". Работа была поддержана грантом РФФИ 18-31-20005 мол-а-вед в § 3, грантом РНФ 18-71-10108 в § 4.}}

{\bf
\copyright\,2019 г.\,\,
М. С. Алкуса$^{*}$,
Д. М. Двинских$^{**, ***}$,
Ф. С. Стонякин$^{*, ****}$,
А. В. Гасников$^{*, **, *****}$,
Д. А. Ковалев$^{******}$
}

*141700 Долгопрудный, М.о., Институтский пер., 9, НИУ МФТИ;

**127051 Москва, Бол. Каретный пер., 19, стр. 1, Ин-т пробл. передачи информац. РАН;

***Институт К. Вейерштрасса, 10117, г. Берлин, Моренштрасса, 39;

****Крымский федеральный университет имени В.И. Вернадского, 295007, г. Симферополь, просп. Академика Вернадского д. 4;

*****Высшая школа экономики, 101000, г. Москва, ул. Мясницкая, 18

******King Abdullah University of Science and Technology, Thuwal, 23955, Saudi Arabia.

e-mail: mohammad.alkousa@phystech.edu,  darina.dvinskikh@wias-berlin.de, fedyor@mail.ru, gasnikov@yandex.ru, dmitry.kovalev@kaust.edu.sa

Поступила в редакцию: 01.12.2019 г.
Исправленный вариант: 27.12.2019 г.
\end{center}

\renewcommand{\abstractname}{\vspace{-\baselineskip}}
\begin{abstract}
В последнее время было показано, как на основе обычного ускоренного градиентного метода решения задач гладкой выпуклой оптимизации можно получить ускоренные методы для более сложных задач (со структурой) и задач, по ходу решения которых используется различная локальная информация о поведении функции (стохастический градиент, гессиан и т.п.). "Ускоренные" \, методы здесь означает, с одной стороны, наличие некоторого единого и достаточно общего способа ускорения. С другой стороны, это означает и "оптимальность" \, методов, что часто удается строго доказать. В настоящей работе предпринята попытка построить в том же духе теорию ускоренных методов  решения гладких выпукло-вогнутых седловых задач со структурой. Основным результатом статьи является получение в некотором смысле необходимых и достаточных условий, при которых сложность решения нелинейных выпукло-вогнутых седловых задач со структурой по числу вычислений градиентов композитов по прямым переменным равна по порядку аналогичной сложности решения билинейных задач со структурой.
\end{abstract}

\textbf{Ключевые слова:} седловая задача, ускоренный метод, слайдинг, проксимально-дружественная функция.

\section{Введение}\label{section:introduction}
Одним из основных направлений в численных методах выпуклой оптимизации в последнее десятилетие стало повсеместное распространение конструкции ускорения обычного градиентного метода, предложенной в 1983 г. Ю.Е. Нестеровым \cite{paper:Nesterov1983}, на различные другие численные методы оптимизации.

Напомним вкратце исходную конструкцию. Для решения задачи выпуклой оптимизации
\begin{equation}\label{problem}
f(x) \rightarrow \min_{x \in \mathbb{R}^{n}},
\end{equation}
было предложено использовать вместо стандартного метода градиентного спуска
$$
    x^{k+1}=x^{k}-\frac{1}{L} \nabla f(x^{k}); \quad  k \geq 0,
$$
где $L$~--- верхняя оценка константы Липшица градиента $f(x)$ в $2$-норме (далее, для краткости, будет говорить, что $f(x)$ имеет $L$--Липшицев градиент), следующий ускоренный метод
$$
	x^{1}=x^{0}-\frac{1}{L} \nabla f(x^{0})
$$
$$
    x^{k+1}=x^{k}-\frac{1}{L} \nabla f\left(x^{k}+\frac{k-1}{k+2}(x^{k}-x^{k-1})\right)+\frac{k-1}{k+2}(x^{k}-x^{k-1}) ; \quad k \geq 1.
$$
Сложность итерации этого метода сопоставима со сложностью итерации градиентного спуска. Однако если градиентный спуск сходится в общем случае (на плохо обусловленных задачах \cite{book:Polyak}) не лучше, чем \cite{paper:Drori_Teboulle}
$$
    f(x^N)-f(x_*) \geq \frac{L R^{2}}{4 N},
$$
то ускоренный метод сходится как
\begin{equation}\label{estimate_AGD}
f(x^N)-f(x_*) \leq \frac{4 L R^{2}}{N^{2}},
\end{equation}
где $x_*$~--- решение задачи \eqref{problem}, а $R = \|x^0 - x_*\|_2$. Если решение не единственное, то под $x_*$ можно понимать то решение, которое наиболее близко (относительно $2$-нормы) к точке старта метода $x^0$ \cite{gasnikov2018book}. Оценка \eqref{estimate_AGD} с точностью до	числового множителя уже не может быть улучшена в общем случае ни на каком другом возможном методе при фиксированном классе задач выпуклой оптимизации \eqref{problem} с $L$--Липшицевым градиентом, см., например,	\cite{Nemirovski:lectures2015,book:Nesterov}. Аналогичные рассуждения можно провести и в невырожденном случае, когда $f(x)$ $\mu$-сильно выпуклая функция \cite{paper:Taylor2017}. В следующих пунктах, мы рассмотрим именно такой случай.

За последние 15 лет описанная конструкция ускорения была успешно перенесена на гладкие задачи условной выпуклой оптимизации, на задачи со структурой (в частности, так называемые, композитные	задачи). Ускорение успешно перенеслось и на неполноградиентные методы (безградиентные методы, спуски по направлению, координатные спуски) и на методы, использующие старшие производные. Также удалось добиться ускорения рандомизированных методов (например, методов редукции дисперсии в задачах минимизации суммы функций) и методов решения задач гладкой стохастической оптимизации. Под успешностью переноса здесь, как и выше, понимается достижение с помощью соответствующих ускоренных методов (с точностью до числовых множителей) известных нижних оценок. Важно также отметить, что в основу схемы ускорения во всех описанных случаях положена идейно одна и та же конструкция. Детали и более подробный обзор литературы можно найти в работе \cite{gasnikov2018book}.

Несмотря на отмеченные достижения, по-прежнему, остается ряд достаточно важных для практики постановок задач, в которых пока не до конца ясно, как именно следует добиваться ускорения имеющихся методов. В частности, можно выделить специальный подкласс задач вида \eqref{problem}, включающий в себя седловые задачи
\begin{equation}\label{eq:3F}
f(x):=r(x)+\underbrace{\max _{y \in Q_{y}}\{F(x, y)-h(y)\}}_{g(x)=F(x, y^*(x))-h(y^*(x))} \rightarrow \min _{x \in Q_{x}},
\end{equation}
где $y^*(x) = \arg\max_{y \in Q_y}\{F(x, y)-h(y)\}$. Заметим, что задачи такого типа встречаются в самых разных приложениях, в том числе в поиске равновесий в двухстадийных моделях транспортных потоков \cite{paper:gasnikov2016searching}.

Если задача не сильно выпукла, то её можно свести к сильно выпуклой применением техники регурялизации и (или) двойственного сглаживания. Без ограничения общности, можно считать $\mu_x \geq \frac{\varepsilon}{2R_x^2}$ (регуляризация \cite{gasnikov2018book}) и (или) $\mu_y \geq \frac{\varepsilon}{2R_y^2}$ (двойственное сглаживание \cite{gasnikov2018book,book:Nesterov,dissertation:Nesterov,book:Lan2019}), где $R_x = \|x^0 - x_*\|_2, R_y = \|y^0 - y^*(x_*)\|_2, (x^0,y^0)$~--- стартовая точка для выбранного численного метода решения задачи \eqref{eq:3F}. Здесь и всюду далее в работе под  $\langle x, y \rangle$ мы понимаем обычное скалярное произведение векторов в конечномерном пространстве и под $ \|x\|_2 = \langle x, x \rangle$ --- евклидову норму.

С точки зрения ускоренных методов данный класс задач достаточно подробно исследован в основном в случае $F(x,y)=\langle Ax,y \rangle$ для некоторого линейного оператора $A$ (см., например, \cite{book:Lan2019}).

В настоящей работе получены аналоги результатов \cite{book:Lan2019} для более общего класса операторов $F$, которые не имеют билинейную структуру. Попутно найдены и некоторые уточнения оценок и для случая $F(x,y)=\langle Ax,y \rangle$. Предлагается сводить рассматриваемые седловые задачи к комбинации вспомогательных задач гладкой сильно выпуклой минимизации отдельно по каждой из групп переменных. Каждая из таких вспомогательных задач может решаться упомянутым выше быстрым градиентным методом Ю. Е. Нестерова. При этом неточность решения вспомогательной (внутренней) подзадачи для одной из групп переменных приводит к необходимости использовать для внешней задачи (которая определяется решением внутренней) концепцию неточного оракула \cite{paper:Devolder2014}, для которой известны оценки скорости сходимости быстрого градиентого метода. В этом случае получены оценки сложности (необходимого количества обращений к подрограмме для вычисления градиента) предложенного метода, аналогичные известному результату об ускоренном градиентном слайдинге Дж. Лана (\cite{book:Lan2019}, раздел 8.2). При этом с одной стороны предложена новая относительно простая схема пояснения ускоренного градиентного слайдинга с помощью техники Каталист \cite{paper:catalyst_2018}, а с другой стороны --- результат об ускоренном градиентном слайдинге обобщён на случай, когда вместо обычных градиентов целевых функционалов используются некоторые их неточные аналоги.

Работа состоит из введения, заключения и трёх основных разделов. В разделе 2 приводится постановка рассматриваемой задачи и кратко описывается подход к ней на базе известного проксимального зеркального метода А.С. Немировского для вариационных неравенств и седловых задач. Раздел 3 посвящён обзору известных результатов о возможности ускорениях оценки скорости сходимости для рассматриваемого класса достаточно гладких сильно выпукло-вогнутых седловых задач, а также формулировкам основных результатов настоящей статьи $(i)$ --- $(iii)$. В разделе 4 приводится схема обоснования основных утверждений работы, сформулированы необходимые вспомогательные утверждения (леммы 1, 2 и 3). Для леммы 3 об ускоренном градиентном слайдинге для минимизации суммы гладких выпуклых функционалов (один из которых сильно выпуклый) при использовании концепции неточного градиента в основной части статьи приведена схема рассмуждений на базе недавно предложенной техники Каталист \cite{paper:catalyst_2018}. Полные доказательства лемм 1, 2 и 3 приводятся в приложении к работе.

\section{Постановка задачи}\label{section:2}

Пусть $Q_x\subseteq \mathbb{R}^m, Q_y\subseteq \mathbb{R}^n$~---непустые, выпуклые и компактные множества, $r: Q_x \to \mathbb{R}$ и $h: Q_y \to \mathbb{R}$~--- $\mu_{x}$-сильно выпуклая и $\mu_y$-сильно выпуклая функции соответственно.

Всюду будем предполагать, что функционал $F: Q_x \times Q_y \to \mathbb{R}$ выпуклый по $x$ и вогнутый по $y$ и задан в некоторой окрестности множества $Q_x \times Q_y$. При этом $F$ будем считать достаточно гладким на $Q_x \times Q_y$. Точнее говоря, для произвольных $x, x' \in Q_x$ и $y, y' \in Q_y$, верны неравенства:
\begin{equation}\label{smooth_F_1}
    \|\nabla_x F(x, y)-\nabla_x F(x', y)\|_2 \leq L_{xx}\|x-x'\|_2,
\end{equation}
\begin{equation}\label{smooth_F_2}
    \|\nabla_x F(x, y)-\nabla_x F(x, y')\|_2 \leq L_{xy}\|y-y'\|_2,
\end{equation}
\begin{equation}\label{smooth_F_3}
    \|\nabla_y F(x, y)-\nabla_y F(x', y)\|_2 \leq L_{xy}\|x-x'\|_2,
\end{equation}
\begin{equation}\label{smooth_F_4}
    \|\nabla_y F(x, y)-\nabla_y F(x, y')\|_2 \leq L_{yy}\|y-y'\|_2,
\end{equation}
где  $L_{xx}, L_{xy}, L_{yy} \geq 0$.

Рассмотрим следующий класс выпукло-вогнутых седловых задач
\begin{equation}\label{problem:min-max}
\min_{x \in Q_{x}}\max_{y \in Q_{y}} \{S(x,y) := r(x) + F(x,y)-h(y)\}.
\end{equation}

Обозначим
\begin{equation}\label{function:S_hat}
    \hat{S}(x,y) = F(x,y) - h(y).
\end{equation}
Тогда можно переписать задачу \eqref{problem:min-max} следующим образом:
$$
    \min_{x \in Q_{x}}\{ r(x) +\max_{y \in Q_{y}} \hat{S}(x,y) \} = \min_{x \in Q_{x}}\{ r(x) + g(x) \},
$$
т.е. задача \eqref{problem:min-max} имеет вид
$$
    f(x):= r(x) + g(x) \to \min_{x \in Q_{x}},
$$
где
\begin{equation}\label{max_problem:g(x)}
    g(x) = \max_{y \in Q_{y}} \hat{S}(x,y).
\end{equation}

Поскольку функционал $\hat{S}(x,\cdot)$ $\mu_y$-сильно вогнутый на $Q_y$, то задача максимизации \eqref{max_problem:g(x)} имеет единственное решение
\begin{equation}\label{solution:max_problem_g(x)}
    y^*(x) := \arg\max_{y \in Q_y} \hat{S}(x,y) \quad
     \forall x \in Q_x,
\end{equation}
откуда
\begin{equation}\label{function:g(x)}
    g(x) = \hat{S}(x,y^*(x)) = F(x,y^*(x)) - h(y^*(x)).
\end{equation}

Всюду далее для произвольного $x \in Q_x$ и некоторого $\delta \geq 0$ будем называть $\tilde{y}_\delta(x) \in Q_y$ $\delta$-приближённым решением задачи \eqref{max_problem:g(x)}, если
\begin{equation}\label{delta_inexact_solution}
    g(x) - \hat{S}(x,\tilde{y}_\delta(x) ) = \hat{S}(x,y^*(x) ) - \hat{S}(x,\tilde{y}_\delta(x) ) \leq \delta.
\end{equation}

Хорошо известно, что задача нахождения седловых точек выпукло-вогнутого функционала может сводиться к задаче решения вариационного неравенства с монотонным оператором
\begin{equation}\label{operator:SPP}
G(x)=\left( \begin{array}{c}{\nabla_{u} f(u, v)} \\ {-\nabla_{v} f(u, v)}\end{array}\right), x=(u, v) \in Q :=Q_{1} \times Q_{2}.
\end{equation}

Напомним общую постановку задачи решения вариационного неравенства (далее сокращённо ВН). Для некоторого оператора $G: Q \to \mathbb{R}^n$, заданного на выпуклом компакте $Q \subset \mathbb{R}^n$ будем рассматривать сильные вариационные неравенства (ВН) вида
\begin{equation}\label{strong_ineq1}
    \left\langle G\left(x_{*}\right), x_{*}-x\right\rangle \leq 0 \quad \forall x \in Q,
\end{equation}
где $G$ удовлетворяет условию Липшица. Отметим, что в \eqref{strong_ineq1} требуется найти $x_* \in Q$ (где $x_*$~--- решение ВН) для которого
\begin{equation}\label{strong_ineq2}
\max _{x \in Q}\left\langle G\left(x_{*}\right), x_{*}-x\right\rangle \leq 0.
\end{equation}

Для решения задачи вариационного неравенства в последние годы весьма популярен проксимальный зеркальный метод А.С. Немировского \cite{paper:Nemirovski2004}. Недавно был предложен также его вариант с адаптивным выбором шага \cite{paper:Gasnikov_Dvurechensky_Stonyakin_Titov_2019}.

Гладкость рассматриваемой седловой задачи приводит к липшицевости оператора $G$, с некоторой константой $L>0$.
 В таком случае, как известно, для проксимального зеркального метода справедлива следующая оценка \cite{paper:Gasnikov_Dvurechensky_Stonyakin_Titov_2019}

\begin{equation}\label{ineq:estimate}
\frac{1}{N} \sum_{k=1}^{N}\left\langle G\left(y^{k}\right), y^{k}-x\right\rangle \leq \frac{L\|x - x^{0}\|_2^2 - L\|x - x^{N}\|_2^2}{2N} \quad \forall x \in Q,
\end{equation}
где
$$
    y^{k} := \arg\min_{x \in Q}\left\{ \langle G(x^{k-1}), x - x^{k-1} \rangle + \frac{L}{2}\|x -x ^{k-1}\|_2^2 \right\}, \quad k = 1,2, \ldots, N.
$$

Заметим, что для всех $y^k \in Q$
\begin{equation}\label{ineq:22}
    \left\langle G\left(x_{*}\right), y^{k}-x_{*}\right\rangle \geq 0.
\end{equation}

Из \eqref{ineq:estimate} следует, что
\begin{equation}\label{ineq:23}
\frac{1}{N} \sum_{k=1}^{N}\left\langle G\left(y^{k}\right), y^{k}-x_{*}\right\rangle \leq \frac{L\left\|x_{*}-x^{0}\right\|_{2}^{2}}{2 N}
\end{equation}

Объединяя неравенства \eqref{ineq:22} и \eqref{ineq:23}, получим
\begin{equation}
\frac{1}{N} \sum_{k=1}^{N}\left\langle G\left(y^{k}\right) - G\left(x_{*}\right), y^{k}-x_{*}\right\rangle \leq \frac{L\left\|x_{*}-x^0\right\|_{2}^{2}}{2 N}
\end{equation}

Если дополнительно предположить сильную монотонность оператора $G$, т.е. существание такого $\mu > 0$, что
\begin{equation}\label{ineq:strongly_monotone}
    \langle G(y) - G(x), y-x\rangle \geq \mu\|y-x\|_{2}^{2} \quad \forall x, y \in Q.
\end{equation}
 С учетом выпуклости функции $\|x\|_2^2$ имеем
\begin{equation*}
    \begin{split}
        \mu\left\|\overline{y}^{N}-x_{*}\right\|_{2}^{2} \leq \mu \sum_{k=1}^{N} \frac{1}{N}\left\|y^{k}-x_{*}\right\|_{2}^{2} & \leq   \frac{1}{N} \sum_{k=1}^{N}\left\langle G \left(y^{k}\right)- G\left(x_{*}\right), y^{k}-x_{*}\right\rangle  \\& \leq \frac{L\left\|x_{*}-x^{0}\right\|_{2}^{2}}{2N},
    \end{split}
\end{equation*}
где $\overline{y}^{N}=\frac{1}{N} \sum_{k=1}^{N} y^{k}$.

На основе последнего неравенства уже возможно организовать процедуру рестартов проксимального зеркального метода и тогда он будет сходится с линейной скоростью. Общая оценка числа итераций, необходимого для достижения приемлемого качества решения будет
\begin{equation}\label{MirrorProxStrong}
\mathrm{O}\left(\frac{L}{\mu} \ln \left(\frac{\mu R^{2}}{\varepsilon}\right)\right),
\end{equation}
где $R = \|x^0 - x_*\|_2$ и $\varepsilon$~--- точность решения $x_*$. Хорошо известно, что данная оценка не может быть улучшена для рассматриваемого класса вариационных неравенств никакими другими методами.

Если применять рассмотренный подход к поставленной задаче \eqref{problem:min-max}, то оператор $G$ из \eqref{operator:SPP} будет $\mu$-сильно монотонным $\mu = \min\{\mu_x, \mu_y\}$, что приведет к такой оценке сложности
\begin{equation}\label{MirrorProxStrongSedlo}
\mathrm{O}\left(\frac{L}{\min\{\mu_x, \mu_y\}} \ln \left(\frac{\min\{\mu_x, \mu_y\} R^{2}}{\varepsilon}\right)\right)
\end{equation}
достижения нужного качества решения. Если $\mu_x$ или $\mu_y$ близко к нулю, то величина в \eqref{MirrorProxStrongSedlo} может оказаться довольно большой достаточно большой. В последующих пунктах мы в частности рассмотрим альтернативные подходы, которые позволяют уточнить оценку \eqref{MirrorProxStrongSedlo} в случае, когда $\mu_x \ll \mu_y$ или $\mu_y \ll \mu_x$.

 \section{Основные результаты}

В этом разделе мы опишем наилучшие известные на данный момент результаты об оценеках скорости сходимости методов для задачи \eqref{eq:3F}, а также сформулируем основные результаты нашей работы.

Будем говорить, что функция $r(x)$ -- проксимально дружественная, если задача вида
\begin{equation}\label{eq:4Fed}
\min _{x \in Q_{x}}\left\{\langle c_{1}, x\rangle + r(x)+ c_{2}\|x\|_{2}^{2}\right\},
\end{equation}
где  $c_1 \in Q_x$  и $c_2 > 0$, может быть решена явно. Аналогично $h(y)$ -- проксимально дружественная функция, если задача вида
\begin{equation}\label{eq:5Fed}
\min_{y \in Q_y}\left\{\langle c_{1}, y\rangle+ h(y)+ c_{2}\|y\|_{2}^{2}\right\},
\end{equation}
может быть решена явно.

Под $\varepsilon$--решением задачи \eqref{eq:3F} будем понимать такую пару $ \left(x^{N_{x}}, y^{N_{y}}\right) \in Q_x \times Q_y$, что
\begin{align*}
f\left(x^{N_x}\right)-f(x_*) \leq & \max_{y \in Q_y \cap B_m(2 R_y)}\left\{r\left(x^{N_x}\right)+F\left(x^{N_x}, y\right)-h(y)\right\}- \\& -\min_{x \in Q_x \cap B_n(2 R_x)}\left\{r(x)+F\left(x, y^{N_y}\right)-h(y^{N_y}) \right\} \leq \varepsilon
\end{align*}
где $B_n(R)$~--- евклидов шар радиуса $R$ в $\mathbb{R}^n$.

Ниже приведены наилучшие известные нам результаты (см. \cite{dissertation:Nesterov,book:Lan2019,Gasnikov_Dvurech_Nesterov:2016,paper:Azizian_2019,paper:Hien_Zhao_Haskell} и цитированную в этих работах литературу) относительно сложности решения задачи \eqref{eq:3F}, которые далее мы постараемся уточнить. Собственно, в пп. 2) -- 4) такое уточнение уже делается: ранее в приведенном виде данные результаты были известны только в предположении $F(x,y)=\langle Ax,y \rangle$ для некоторого линейного оператора $A$. В этом случае $L_{xx} = L_{yy} = 0$, $L_{xy} = L_{yx} = \sqrt{\lambda_{\max}(A^T A)}$, где $\lambda_{\max}(A^T A)$~--- наибольшее собственное значение матрицы $A^TA$.

\begin{enumerate}
	\item[1)] Если $r(x)$ и  $h(y)$ проксимально дружественны, то $\varepsilon$--решение задачи \eqref{eq:3F} может быть достигнуто, за  $\tilde{O}\left(\frac{L_{xy}}{\sqrt{\mu_x \mu_y}} \right)$ вычислений \eqref{eq:4Fed}, $\nabla_x F(x,y)$  и \eqref{eq:5Fed}, $\nabla_y F(x,y)$ при $F(x,y)=\langle Ax,y \rangle$ и за $\tilde{O}\left(\frac{\max \{L_{xx},L_{xy},L_{yy}\}}{\min\{\mu_x, \mu_y\}} \right)$ вычислений \eqref{eq:4Fed}, $\nabla_x F(x,y)$  и \eqref{eq:5Fed}, $\nabla_y F(x,y)$ в общем случае. Здесь и далее $\tilde{O}() = O()$ с точностью до небольшой степени логарифмического по $\varepsilon, \mu_x$ или $\mu_y$, а также по $R_x$ или $R_y$ множителя. Как правило, показатель этой степени 1 или 2.
	\item[2)] Если $r(x)$ имеет $L_x$--Липшицев градиент, но не проксимально дружественна, то $\varepsilon$--решение задачи \eqref{eq:3F} может быть достигнуто за $\tilde{O}\left(\sqrt{\frac{L_x}{\mu_x}}\right)$ вычислений  $\nabla r(x)$,
	\begin{equation}\label{est:1}
	\tilde{O}\left(\sqrt{\frac{L_{xx}+\frac{L_{xy}^2}{\mu_y}}{\mu_x}}\right)
	\end{equation}
	вычислений $\nabla_x F(x,y)$ и
	\begin{equation}\label{est:2}
	\tilde{O}\left(\sqrt{\frac{L_{xx}+\frac{L_{xy}^2}{\mu_y}}{\mu_x}} \sqrt{\max \left\{ \frac{L_{yy}}{\mu_y}, 1 \right\}  }   \right)
	\end{equation}
	вычислений \eqref{eq:5Fed}, $\nabla_y F(x,y)$.
	\item[3)] Если $h(y)$ имеет $L_y$--Липшицев градиент, но не проксимально дружественна, то $\varepsilon$--решение задачи \eqref{eq:3F} может быть достигнуто за
	\eqref{est:1} вычислений \eqref{eq:4Fed}, $\nabla_xF(x,y)$ и
	\begin{equation}\label{est:3}
	\tilde{O}\left(\sqrt{\frac{L_{xx}+\frac{L_{xy}^2}{\mu_y}}{\mu_x}} \sqrt{\frac{L_{y}}{\mu_y} } \right)
	\end{equation}
	вычислений $\nabla h(y)$, \eqref{est:2}
	вычислений $\nabla_y F(x,y)$.
	\item[4)] Если $r(x)$ имеет $L_x$--Липшицев градиент, $h(y)$ имеет $L_y$--Липшицев градиент, но обе функции не проксимально дружественны, то $\varepsilon$--решение задачи \eqref{eq:3F} может быть достигнуто за  $\tilde{O}\left(\sqrt{\frac{L_x}{\mu_x}}\right)$ вычислений  $\nabla r(x)$,
	\eqref{est:1} вычислений $\nabla_x F(x,y)$ и
	\eqref{est:3} вычислений $\nabla h(y)$,
	\eqref{est:2} вычислений $\nabla_y F(x,y)$.
\end{enumerate}

Отметим, что результаты всего п. 1) и пп. 2), 4) в части числа вычислений $\nabla r(x)$ (и, по-видимому, в части числа вычислений $\nabla_x F(x,y)$) в общем случае не могут быть улучшены \cite{Nemirovski:lectures2015,paper:Ouyang2018} (с точностью до логарифмического множителя), если $$\dim(x)+\dim(y) \gg \frac{L_{xy}}{\sqrt{\mu_x \mu_y}}$$
Заметим, что правую часть можно менять в зависимости от специфики постановки задачи. Если это условие не выполняется, то можно свести седловую задачу к негладкой задаче выпуклой оптимизации \cite{book:Gadan_3,book:Nemirovski_Ydin_1979} и решать ее методами типа центров тяжести, например, методом эллипсоидов или Вайды \cite{book:Nemirovski_Ydin_1979,boobk:Bubeck_2015}. Заметим, что методом эллипсоидов можно решать седловую задачу и напрямую \cite{paper:Nemirovski_Onn_2010}.

Сформулируем результаты настоящей работы.

\begin{itemize}
    \item[i)]
	Основным результатом работы является обоснование возможности в пп. 2)--4) выше полностью убрать оговорку "при $F(x,y) = \langle Ax,y \rangle$".
	
	При $\mu_x = \frac{\varepsilon}{R_x^2}$ частично это уже было сделано в работе \cite{paper:Hien_Zhao_Haskell}. Однако полученные в этой работе оценки в части числа вычислений  $\nabla r(x)$ проигрывают оценкам, приведенным выше в $\sim\frac{1}{\mu_y}$ раз.
	
	\item[ii)]
	Другим результатом работы является уточнение приведенных выше утверждений в случае, когда $F(x,y) = \langle Ax,y \rangle$ для некоторого линейного оператора $A$, $Q_y = \mathbb{R}^m,  h(y)$ имеет $L_y$--Липшицев градиент и $\frac{\lambda_{\min }\left(A^{T} A\right)}{L_{y}} \gg \mu_{x}$. В этом случае во всех приведенных выше формулах можно заменить $\mu_x$ на $\frac{\lambda_{\min}\left(A^{T} A\right)}{L_{y}}$.
	
	\item[iii)]
	Более того, если $F(x,y) = \langle Ax,y \rangle$ для некоторого линейного оператора $A$, $Q_x = \mathbb{R}^n, Q_y = \mathbb{R}^m$, $r(x)$ имеет $L_x$--Липшицев градиент, $h(y)$ имеет $L_y$--Липшицев градиент, то оценки на количество вычислений $\nabla_{x} F(x,y)= A^T y, \nabla_{y} F(x,y)= Ax$ и \eqref{eq:5Fed}  (здесь предполагается проксимальная дружественность $h(y)$) можно заменить (это имеет смысл, если новые оценки станут лучше имеющихся) на следующую оценку
	\begin{equation}\label{estimateiii}
	    [\text{Число вычислений} \; \nabla r(x)]\cdot \tilde{O}\left( \sqrt{\frac{L_y \lambda_{\max}(A^TA)}{\mu_y \lambda_{\min }^+(A^TA)}} \right) = \tilde{O}\left( \sqrt{\frac{L_x L_y \lambda_{\max}(A^TA)}{\mu_x \mu_y \lambda_{\min }^+(A^TA)}} \right),
	\end{equation}
	\end{itemize}
где $\lambda_{\min}^+(A^TA)$--минимальное положительное собственное значение матрицы $A^TA$.

Отметим, что аналогичное ii) утверждение можно сформулировать и в случае, когда $r(x)$ имеет $L_x$--Липшицев градиент.

Стоит отметить, что результат ii) позволяет, среди прочего, элементарным образом объяснить, как решать матричную игру за $\tilde{O}\left( \sqrt{\frac{ \lambda_{\max }(A^TA)}{\lambda_{\min}(A^TA)}} \right)$ матрично-векторных умножений. Сначала необходимо двойственно сгладить (должным образом) исходную постановку задачи \cite{gasnikov2018book,book:Nesterov,dissertation:Nesterov,book:Lan2019}, т.е. ввести проксимально-дружественную функцию $h(y) = \frac{\varepsilon \|y\|_2^2}{4R_y^2}$. Отметим, что при этом $Q_x = Q_y =\mathbb{R}^n$. Далее, полученную в результате функцию $g(x)$ можно минимизировать обычным ускоренным (быстрым градиентным) методом \cite{paper:Nesterov1983,gasnikov2018book,Nemirovski:lectures2015,book:Nesterov,paper:Taylor2017,dissertation:Nesterov,book:Lan2019} (в сильно выпуклом случае). Полученная при этом сложность решения задачи существенно лучше популярных и активно исследуемых процедур типа экстраградиентного метода, для которого удалось лишь получить такую оценку трудоемкости $\tilde{O}\left( \frac{ \lambda_{\max }(A^TA)}{\lambda_{\min}(A^TA)} \right)$ \cite{paper:Mokhtari2019}.

Отметим также, что результаты, приведенные в п. iii) нашли важное приложение в построении оптимальных ускоренных методов для гладких выпуклых задач децентрализованной распределенной оптимизации \cite{paper:Dvinskikh_Gasnikov2019}. Собственно, наш изначальный интерес к этой проблематике и был связан с развитием идей работы \cite{paper:Dvinskikh_Gasnikov2019}.


\section{Необходимые вспомогательные утверждения и схема обоснования основных результатов работы}\label{section:3}

Для того чтобы получить отмеченные в предыдущем пункте результаты, приведем необходимые и, в основной части, уже известные результаты.	

\begin{lemma}[\cite{gasnikov2018book,paper:Kakde_Shalev_Tewari,book:Rockafellar}]\label{lemma1}
	В обозначениях п. 2, если $F(x,y) = \langle Ax,y \rangle$, то $g(x)$ будет иметь $L$--Липшицев градиент, где $L = \frac{\lambda_{\max }(A^TA)}{\mu_y}$. Если, дополнительно, $h(y)$ имеет $L_y$--Липшицев градиент и $Q_y = \mathbb{R}^m$, то $g(x)$~--- $\mu$-сильно выпуклая функция на $(\text {Ker} A)^{\perp}$ где $\mu = \frac{\lambda_{\min }^+(A^TA)}{L_y}$ При этом $\nabla g(x) \in (\text {Ker} A)^{\perp}$.
\end{lemma}

Следующее утверждение демонстрирует необходимость использования концепции неточного градиента для рассматриваемого подхода к седловым задачам.

\begin{lemma}[\cite{paper:Devolder2014,Gasnikov_Dvurech_Nesterov:2016,paper:Hien_Zhao_Haskell}] \label{lemma2}
	В обозначениях п. 2, $g(x)$ будет иметь $L$--Липшицев градиент $\nabla g(x) = \nabla_{x}F(x,y^*(x))$ (теорема Демьянова--Данскина), где $L = L_{xx} + \frac{2L_{xy}^2}{\mu_y}$. Более того, для всех $x, z \in Q_x$ выполняется следующее условие (словами: $\nabla_x F(x, \tilde{y}_\delta(x))$~--- $(2\delta,2L)$--градиент для $g(x)$)
	\begin{equation}\label{eq:lemma2}
	0 \leq g(z) -  \left[\left\{ F(x, \tilde{y}_{\delta}(x))-h(\tilde{y}_\delta(x)) \right\}
	+\langle \nabla_{x}  F(x, \tilde{y}_\delta(x)), z-x \rangle\right]
	\leq \frac{2 L}{2}\|z-x\|_{2}^{2}+2 \delta,
	\end{equation}
	где $\tilde{y}_\delta(x)$ такой, что
	$$
	    \underbrace{\max_{y \in Q_{y}} \{F(x, y)-h(y)\}}_{g(x)} - \left\{F\left(x,\tilde{y}_{\delta}(x)\right)-h\left(\tilde{y}_{\delta}(x)\right)\right\} \leq \delta.
	$$
	который определен в \eqref{delta_inexact_solution}.
\end{lemma}

Рассмотрим наиболее тонкое из используемых нами вспомогательных утверждений \cite{book:Lan2019}. Пусть нужно решить задачу сильно выпуклой композитной оптимизации
	\begin{equation}\label{problem_lemma3}
	r(x)+g(x) \rightarrow \min_{x\in \mathbb{R}^{n}},
	\end{equation}
	c точностью $\varepsilon$ по функции. Считаем, что функции $r(x)$ и $g(x)$ имеют константы Липшица градиента $L_r$ и $L_g$, и хотя бы одна из этих функций $\mu$-сильно выпуклая. К такой задаче  \eqref{problem_lemma3} можно применить результат о так называемом ускоренном градиентном слайдинге из п. 8.2 \cite{book:Lan2019}.
	\begin{lemma} \label{lemma3}
	Необходимого качества решения задачи решения \eqref{problem_lemma3} можно достичь за $N_r = \tilde{O} \left(\sqrt{\frac{L_r}{\mu}}\right)$ вычислений градиента первого слагаемого $\nabla r(x)$ и $N_g = \tilde{O}\left(\sqrt{\frac{L_g}{\mu}}\right)$ вычислений градиента второго слагаемого $\nabla g(x)$. Результат останется верным, если вместо $\nabla g(x)$ использовать $\left(O\left(\frac{\varepsilon}{N_g}\right), O(L_g)\right)$-градиент.
\end{lemma}

Подробное доказательство леммы \ref{lemma3} приведено в приложении. Здесь мы лишь приведем некоторую схему рассуждений. Для наглядности полагаем $\mu$-сильно выпуклым именно функционал $g$. Если это не так и $r$~--- $\mu_r$-сильно выпуклый, а $g$ просто выпуклый, то можно заменить $r(x)$ на $r(x) - \frac{\mu_g}{2} \|x\|_2^2$ и $g(x)$ на $g(x) + \frac{\mu_g}{2} \|x\|_2^2$ для некоторого $\mu_g < \mu_r$. Тогда оба функционала будут как гладкими, так и сильно выпуклыми.

К рассмотренной задаче \eqref{problem_lemma3} можно применить технику Каталист (алгоритм \ref{catalyst}) \cite{paper:catalyst_2018}, в предположении $L_r \ll L_g$ и $\mu$-сильной выпуклости $g$ в $2$-норме, причем $\mu \ll L_r$.

\begin{algorithm}[h]
	\caption{Каталист}
	\begin{algorithmic}[1]
		\STATE {\bf Вход:} $x^0 \in \mathbb{R}^n,$ параметр $L$.
		\STATE $y^0 := x^0$
		\WHILE{желаемая точность не достигнута}
			\STATE Найти $x^k$ с некоторой точностью применением алгоритма \ref{alg:gradient_for_composite}
			$$
			    x^{k} \approx \underset{x \in \mathbb{R}^n}{\arg\min }\left\{r(x)+g(x)+\frac{L}{2}\|x-y^{k-1}\|_2^{2}\right\}
			$$
			\STATE Вычисляем $y^k$ используя шаг экстраполяции, с $\beta_k \in (0,1)$
			$$
			    y^k = x^k + \beta_k(x^k - x^{k-1}).
			$$
		\ENDWHILE
		\STATE {\bf Выход:} $x^k$.
	\end{algorithmic}
	\label{catalyst}
\end{algorithm}
\begin{algorithm}[h]
	\caption{Неускоренный градиентный метод для задач композитной оптимизации.}
	\begin{algorithmic}[1]
		\STATE {\bf Вход:} $x^0 \in \mathbb{R}^n,$ параметр $L_r$.
		\FOR{$k = 0, 1, 2, \ldots$}
			\STATE
			$$
			    \phi_{k+1}(x):= \langle \nabla r(x^k), x - x^k \rangle + g(x) - g(x^k) + \frac{L_r}{2} \|x - x^k\|_2^2,
			$$
			$$
			   x^{k+1}:=\arg\min_{x \in \mathbb{R}^n}\phi_{k+1}(x)
			$$
		\ENDFOR
		\STATE {\bf Выход:} $\bar{x}_N=\frac{1}{N} \sum_{k=0}^{N-1} x^{k+1}$.
	\end{algorithmic}
	\label{alg:gradient_for_composite}
\end{algorithm}

Тогда (см. \cite{paper:catalyst_2018}) вместо исходной задачи потребуется $\tilde{O}\left(\sqrt{\frac{L}{\mu}}\right)$ раз решать задачу вида
\begin{equation}\label{AuxProblem}
r(x)+g(x)+\frac{L}{2}\left\|x-y^{k-1}\right\|_{2}^{2} \rightarrow \min _{x\in \mathbb{R}^{n}},
\end{equation}
где $L > 0$~--- некоторый параметр регуляризации, а последовательность $y^0, y^1, \ldots $ образуется согласно схеме алгоритма \ref{catalyst}. При этом задачу \eqref{AuxProblem} можно решать неускоренным композитным градиентным методом (алгоритм \ref{alg:gradient_for_composite}), считая $g(x)+\frac{L}{2}\left\|x-y^{k-1}\right\|_{2}^{2}$ композитом. Как известно, число итераций такого метода будет совпадать с количеством вычислений $\nabla r(x)$ и равно $\tilde{O}\left( \frac{L_r}{L + \mu}  \right)$. Но при этом не предполагается проксимальная дружественность функции $g(x)$ и поэтому необходимо учитывать сложность решения возникающей на каждой итерации неускоренного композитного градиентного метода задачи вида
$$
\left\langle\nabla r(\tilde{x}^{l}), x-\tilde{x}^{l}\right\rangle+\frac{L_r}{2}\left\|x-\tilde{x}^l\right\|_2^2+g(x)+\frac{L}{2}\left\|x-y^{k-1}\right\|_2^2 \rightarrow \min_{x \in \mathbb{R}^n},
$$
Для решения уже этой вспомогательной задачи можно использовать ускоренный композитный градиентный метод для задач сильно выпуклой оптимизации \cite{paper:Gasnikov_MIPT_2016} (см. также алгоритмы 3 -- 4 далее). При этом слагаемое $\frac{L_r}{2}\left\|x-\tilde{x}^l\right\|_2^2+\frac{L}{2}\left\|x-y^{k-1}\right\|_2^2$ считается композитом. Число итераций такого метода будет $\tilde{O}\left(\sqrt{\frac{L_g}{L_r + L + \mu} }\right)$. Таким образом, общее количество вычислений $\nabla g(x)$ будет
$$
\tilde{O}\left(\sqrt{\frac{L}{\mu}} \right) \cdot \tilde{O}\left( \frac{L_r}{L + \mu}\right) \cdot \tilde{O}\left(\sqrt{\frac{L_g}{L_r + L + \mu}}\right).
$$
Выберем параметр регуляризации $L$ так, чтобы последнее выражение было минимальным. Тогда с учетом сделанных предположений $L_r \ll L_g$ и  $\mu \ll L_r$ получим, что $L \simeq L_r$. Следовательно, общее число вычислений $\nabla g(x)$ будет действительно равно $\tilde{O}\left(\sqrt{\frac{L_g}{\mu}} \right)$.

Заметим, что полностью новым является утверждение из последнего предложения формулировки леммы \ref{lemma3}. Доказательство этого ключевого наблюдения проводится аналогично оригинальной статье \cite{paper:Devolder2014}, см. также \cite{book:Lan2019,paper:Hien_Zhao_Haskell}. Действительно,
если функционал $g:\mathbb{R}^n \rightarrow \mathbb{R}$ допускает
($\delta, L$)-градиент $\nabla_{\delta} g(x)$ в любой запрошенной точке $x$, т.е. верно неравенство
$$
    g(y) \leq g(x)+\left\langle \nabla_{\delta}g(x), y - x \right\rangle + \frac{L}{2}\|y-x\|_2^2+\delta.
$$
Последнее неравенство отличается от стандартного условия $L$-липшицевости градиента выпуклого функционала $g$
\begin{equation}\label{LipGrad}
g(y) \leq g(x)+\left\langle \nabla g(x), y - x \right\rangle + \frac{L}{2}\|y-x\|_2^2
\end{equation}
лишь постоянной величиной погрешности $\delta$. Поэтому применение указанного неравенства для модификации результата об ускоренном градиентном слайдинге \cite{book:Lan2019} (который в стандартном случае оснван на комбинации оценок, связанных с неравенством \eqref{LipGrad}) для решения $N_g$ вспомогательных подзадач приведет к накоплению в итоговой оценке погрешности, сопоставимой с величиной $O(N_g \delta)$. Это обстоятельство и приводит к необходимости вместо $\nabla g(x)$ использовать имено $\left(O\left(\frac{\varepsilon}{N_g}\right), O(L_g) \right)$-градиент.

Утверждения i) и ii) из предыдущего раздела работы получаются сочетанием лемм \ref{lemma1}--\ref{lemma3}. Приведём схему доказательства первого из результатов. Напомним, что мы рассматриваем задачу вида
\begin{equation}\label{minimax_problem}
	\min_{x \in Q_{x}}\max_{y \in Q_{y}} \{S(x,y) = r(x) + F(x,y) - h(y)\}.
\end{equation}
Мы проводим доказательство в предположении, что функционалы $r$ и $h$ проксимально дружественны (простой структуры, т.е. композиты).
\begin{remark}
Если некоторый функционал $S(x,y)$ $\mu_x$-сильно выпуклый по $x$ и $\mu_y$-сильно вогнутый по $y$, то можно представить $S$ в виде
$$
    S(x,y) = F(x,y) + \frac{\mu_x}{2} \|x\|_2^{2} - \frac{\mu_y}{2} \|y\|_2^{2},
$$
причём функционал $F$ сохранит свойство выпуклости по $x$ и вогнутости по $y$.
Отметим, что в виду достаточной гладкости функционалов $\|x\|_2^{2}$ и $\|y\|_2^{2}$ (липшицевость градиентов) если $S$ гладкий, то $F$ будет гладким. Таким образом, приводимая ниже схема рассуждений применима к седловым задачам для произвольного $\mu_x$-сильно выпуклого по $x$ и $\mu_y$-сильно вогнутого по $y$, а также достаточно гладкого функцонала $S(x, y)$.
\end{remark}

Зафиксируем $x \in Q_x$ и введём вспомогательный функционал
\begin{equation}\label{eq:4F}
f(x) := \max\limits_{y \in Q_y} S(x,y),
\end{equation}
что позволяет рассматривать седловую задачу \eqref{minimax_problem} как задачу выпуклой минимизации
\begin{equation}\label{eq:5F}
f(x) \to \min\limits_{x \in Q_x}.
\end{equation}

Ясно, что
$$
    f(x) = \max\limits_{y \in Q_y}\{ r(x) + F(x,y) - h(y) \} = r(x) + \underbrace{\max _{y \in Q_{y}}\{F(x, y)-h(y)\}}_{g(x)=F(x, y^*(x))-h(y^*(x))},
$$

По лемме \ref{lemma2} функционал $f$ (см. \eqref{eq:4F}) $\mu_x$-сильно выпуклый и имеет $L$--Липщицев градиент
\begin{equation}\label{eq:6}
\| \nabla f(x_2) - \nabla f(x_1)\|_2   \leq L \|x_2 - x_1\|_2, \quad \forall x_1, x_2 \in Q_x,
\end{equation}
где
\begin{equation}\label{eq:8}
L = L_{xx} + \frac{2 L_{xy}^2}{\mu_y}.
\end{equation}
Также согласно лемме \ref{lemma2} (см. также \cite{paper:Hien_Zhao_Haskell}) для функционала $f$ выполнено неравенство:
\begin{equation}\label{eq:7}
f(x_2) \leq S(x_1, \tilde{y}_{\gamma}(x_1) ) + \langle \nabla_{x} S(x_1, \tilde{y}_{\gamma}(x_1) ), x_2  - x_1 \rangle   + L \|x_2 - x_1\|_2^2 + 2 \gamma,
\end{equation}
а $\gamma$~--- точность решения вспомогательной задачи максимизации $S(x_1,y)$ по $y$:
\begin{equation}\label{eq:9}
f(x_1) - S(x_1, \tilde{y}_{\gamma}(x_1) ) \leq \gamma.
\end{equation}
Задача $\mu_y$--сильно вогнутой композитной (в силу предположения для $h$ выше) максимизации
$$
    F(x,y)-h(y) \to \max\limits_{y \in Q_y}
$$
(при фиксированном $x \in Q_x$) может быть решена с точностью  $\gamma$ в \eqref{eq:9}  за
\begin{equation}\label{eq:10}
O\left(\sqrt{\frac{L_{yy}}{\mu_y}} \ln \frac{1}{\gamma} \right)
\end{equation}
итераций быстрого градиентного метода.

Покажем, как с использованием рестартов быстрого градиентного метода в концепции $(\delta, L)$-оракула можно достичь заданного качества решения по функции для задач композитной оптимизации с требуемой оценкой сложности. Для задачи композитной оптимизации
\begin{equation}\label{problem:composite_opt}
    f(x):= r(x)+g(x) \rightarrow \min_{x\in Q}.
\end{equation}
можно применить предложенный в \cite{paper:gasnikov_turin_2019} алгоритм \ref{alg:acculerated_gradient_for_composite}.

\begin{algorithm}[t]
	\caption{Быстрый градиентный метод с $(\delta, L)$-оракулом для задач композитной оптимизации \eqref{problem:composite_opt}.}
	\begin{algorithmic}[1]
		\STATE {\bf Вход:} $x^0 \in Q$~--- начальная точка, $N$~--- количество шагов, $L> 0$ и $\delta > 0$.
		\STATE $y^{0}:=x^{0}, \quad u^{0}:=x^{0}, \quad \alpha_{0}:=0, \quad A_{0}:=0$
		\FOR{$k = 1, 2, \ldots, N$}
			\STATE Находим наибольший корень, $\alpha_{k+1}$ так что
			$$ A_k + \alpha_{k+1} = L\alpha^2_{k+1},$$
			\STATE $$A_{k+1} := A_k + \alpha_{k+1},$$
			\STATE $$y^{k+1} := \frac{\alpha_{k+1}u^k + A_k x^k}{A_{k+1}}, $$
			\STATE \begin{gather*}
                    \phi_{k+1}(x) = \frac{1}{2}\|x- u^k\|_2^2 + \alpha_{k+1}\left( \langle \nabla r(y^{k+1}), x- y^{k+1} \rangle + g(x) - g(y^{k+1}) \right), \\
                    u^{k+1} := {\arg\min_{x \in Q}} \,\phi_{k+1}(x),
                  \end{gather*}
			\STATE $$x^{k+1} := \frac{\alpha_{k+1}u^{k+1} + A_k x^k}{A_{k+1}},$$
		\ENDFOR
		\STATE {\bf Выход:} $x^N$.
	\end{algorithmic}
	\label{alg:acculerated_gradient_for_composite}
\end{algorithm}

\begin{algorithm}[h]
	\caption{Быстрый градиентный метод для задач сильно выпуклой композитной оптимизации с $(\delta, L)$-оракулом, рестарты алгоритма \ref{alg:acculerated_gradient_for_composite}.}
	\begin{algorithmic}[1]
		\STATE {\bf Вход:} $x^0 \in Q$~--- начальная точка, $L>0, \mu>0, \delta>0,$ $\varepsilon$~--- точность решения, $R$  и  $p= \left\lceil \log_2\left( \frac{\mu R^2}{\varepsilon}\right) \right\rceil$~--- количество рестартов.
		\FOR{$j = 1, \ldots, p$}
			\STATE Выполнить $N_j = \left\lceil 3 \sqrt{\frac{2L}{\mu}} \right\rceil$ итераций алгоритма \ref{alg:acculerated_gradient_for_composite} .
			\STATE $x^0 := x^{N_j}$.
		\ENDFOR
		\STATE {\bf Выход:} $\hat{x}:=x^{N_p}$.
	\end{algorithmic}
	\label{alg:restart}
\end{algorithm}

Для алгоритма \ref{alg:acculerated_gradient_for_composite} в \cite{paper:gasnikov_turin_2019} доказана следующая оценкка скорости сходимости (мы предполагаем, что вспомогательные задачи решаются точно):
$$
    f(x^N)-f(x_*) \leq \frac{8 L R^2}{(N+1)^2} + 2 N \delta,
$$
где $R^2 = \frac{1}{2}\|x_*- x^0\|_2^2$, а $x_*$~--- ближайшая точка минимума к точке $x^0$. Введём вспомогательное обозначение
\begin{equation}
     \psi(x,y):= \langle\nabla r(y), x - y \rangle +g(x) - g(y).
\end{equation}

Тогда ввиду $\mu_x$-cильной выпуклости $g$ имеем
\begin{equation}
\frac{\mu_x}{2} \|x-y\|_2^2 \leq f(x)-\left(f(y) + \psi(x,y)\right) \leq \frac{L}{2} \|y-x\|_2^2 + \delta.
\end{equation}

Если сделать естественное предположение $\psi(x, x_*) \geq \langle \nabla(r+g)(x_*), x - x_* \rangle\geq 0$, то после $N_1$ итераций алгоритма \ref{alg:acculerated_gradient_for_composite} имеем
\begin{equation}
    \frac{\mu_x}{2}\|x^{N_1} - x_* \|_2^2 \leq f\left(x^{N_1}\right)-f\left(x_*\right) \leq \frac{4L \|x^0 - x_* \|_2^2}{N_{1}^2}+2 N_1 \delta.
\end{equation}
Если для числа итераций $N_1$ алгоритма \ref{alg:acculerated_gradient_for_composite} выбрать $\delta$ так, чтобы выполнялось неравенство
\begin{equation}
    2 N_1 \delta \leq  \frac{L \|x^0 - x_* \|_2^2}{2 N_{1}^2},
\end{equation}
то получим
\begin{equation}
     \|x^{N_1} - x_*\|_2^2 \leq  \frac{9L}{\mu_x N_1^2}\|x^0 - x_*\|_2^2.
\end{equation}
Поэтому, выбирая  $N_1 = \left\lceil 3e \sqrt{\frac{L}{\mu_x}}\right\rceil$\footnote{Данный способ выбора количества итераций авторам подсказал Роланд Хильдебранд.}, получим
$$
    \|x^{N_1} - x_*\|_2^2 \leq \frac{1}{e^2}\|x^0 - x_*\|_2^2.
$$

После этого выберем для aлгоритма \ref{alg:acculerated_gradient_for_composite} в качестве точки старта $x^{N_1}$, и снова сделаем $N_1$ итераций, и т.д. Ясно, что для достижения приемлемого качества решения можно выбрать количество рестартов $p$ алгоритма \ref{alg:acculerated_gradient_for_composite} для алгоритма \ref{alg:restart} следующим образом:
$$
    p = \left\lceil \frac{1}{2} \ln\left( \frac{\mu_x R^2}{\varepsilon}\right) \right\rceil,
$$
где $R^2 = \frac{1}{2} \|x^0 - x_*\|_2^2$. В таком случае общее число итераций алгоритма \ref{alg:restart} в предложенной схеме будет
$$
    N = \left\lceil \ln\left( \frac{\mu_x R^2}{\varepsilon}\right) \right\rceil \cdot \left\lceil \frac{3e}{2} \sqrt{\frac{L}{\mu_x}} \right\rceil,
$$
т.е.
$$
    N = O\left(\sqrt{\frac{L}{\mu_x}} \left\lceil \ln\left( \frac{\mu_x R^2}{\varepsilon}\right) \right\rceil \right).
$$

Итак, ясно, что при $\delta=O(\varepsilon)$ предложенной схемой возможно получить точность $f(\hat{x})-f(x_*)\leq\varepsilon$ для выхода $\hat{x}$ алгоритма \ref{alg:restart}.

Далее, неравенство \eqref{eq:7} (см. также неравенство (2.20) из \cite{paper:Hien_Zhao_Haskell}) означает, что функционал $f$ в произвольной точке допускает $(2\gamma, L, \mu_x)$-оракул в смысле \cite{thesis:devolder2013}.

Таким образом, точность $\gamma$ для решения внутренней задачи (аналог точности $\delta$ в концепции ($\delta, L$)-градиента) нужно выбирать как $O(\varepsilon)$ и тогда для внешней задачи будет гарантирована точность $\varepsilon$ по функции.
Это означает, что общее количество вызовов оракула для $\nabla_x S(\cdot, \cdot)$ при решении внешней задачи \eqref{eq:5F} будет равно
\begin{equation}\label{eq:14}
O\left(\sqrt{\frac{L_{xx}}{\mu_x } + \frac{2L_{xy}^2}{\mu_x \mu_y}} \ln \frac{1}{\varepsilon}\right).
\end{equation}
Количество же вызовов оракула для $\nabla_y S(\cdot, \cdot)$ при решении внутренней задачи будет равно в силу \eqref{eq:10}
\begin{equation}\label{eq:13}
O \left(\sqrt{\frac{L_{yy}}{\mu_y}}\sqrt{\frac{L_{xx}}{\mu_x} + \frac{2L_{xy}^2}{\mu_x \mu_y} } \ln^2 \frac{1}{\varepsilon}\right).
\end{equation}

\begin{remark}
	Согласно известной теореме о минимаксе сильная выпукло-вогнутость функционала
	$S(\cdot, \cdot)$ означает, что:
	$$
		\min\limits_{x \in Q_x}\{ f(x) = \max\limits_{y \in Q_y} S(x,y)\} = \max\limits_{y \in Q_y}\{ l(x) = \min\limits_{x \in Q_x} S(x,y)\}
	$$
	и можно свести задачу \eqref{minimax_problem} к задаче вогнутой максмимзации
	\begin{equation}\label{eq:15}
	\tilde{f}(y) := \min\limits_{x \in Q_x}S(x,y) \to \max\limits_{y \in Q_y}.
	\end{equation}
	В таком случае количество вызовов $\nabla_y S(\cdot, \cdot)$ будет
	\begin{equation}\label{eq:16}
	O \left(\sqrt{\frac{L_{yy}}{\mu_y} + \frac{2L_{xy}^2}{\mu_x \mu_y} } \ln \frac{1}{\varepsilon}\right),
	\end{equation}
	а количество вызовов $\nabla_xS(\cdot, \cdot)$
	$$
		O \left(\sqrt{\frac{L_{xx}}{\mu_x}}\sqrt{\frac{L_{yy}}{\mu_x \mu_y} + \frac{2L_{xy}^2}{\mu_x^2 \mu_y} } \ln^2 \frac{1}{\varepsilon}\right).
	$$
\end{remark}

Утверждение i) в случае, когда какая-то из функций $r(x)$ или $h(y)$ не проксимально дружественна, получается похожей схемой рассуждений. Но при этом мы получим вспомогательные подзадачи вида \eqref{problem_lemma3} с сильно выпуклой целевой функцией, которые можно решать с использованием результата об ускоренном градиентном слайдинге. Сформулированные оценки утверждения i), таким образом, вытекают из леммы \ref{lemma3}.

Отметим, что для утверждения леммы \ref{lemma3} о сложности решения вспомогательной подзадачи \eqref{problem_lemma3} в случае сильной выпуклости $r$ и $g$ не имеет значения, параметр сильной выпуклости какого из этих функционалов использовать в оценках из леммы \ref{lemma3}. Поэтому утверждение ii) теперь будет следовать из леммы \ref{lemma1} (где приводится оценка на параметр сильной выпуклости $g$).

Скажем несколько слов про обоснование утверждения iii). Мы отправляемся от \cite{paper:Dvinskikh_Gasnikov2019}.

Ясно, что решаемую задачу минимизаци функции $f(x)= r(x) + g(x)$ на множестве $Q_x = \mathbb{R}^n$ можно рассматривать как задачу композитной оптимизации с композитом $g$. Точнее говоря, ввиду $\mu_x$-сильной выпуклости и $L_x$--липшицевости градиента $r$ можно свести задачу минимизации $f$ к рассмотрению семейства $k = \tilde{O}\left(\sqrt{\frac{L_x}{\mu_x}}\right)$ вспомогательных задач вида
\begin{equation}\label{statiiiaux}
\alpha_i (\langle \nabla r(x_i), x - x_i \rangle + g(x) - g(x_i)) + \frac{L_x \|x-x_i\|_2^2}{2} \rightarrow \min\limits_{x \in Q_x}, \quad i = 1, 2, ..., k.
\end{equation}
Поскольку $g$ уже, вообще говоря, не имеет простой структуры, то необходимо учитывать сложность каждой из $k$ задач \eqref{statiiiaux}. Поскольку вдоль ядра $ Ker A$ функия $g$ принимает постоянное значение, то без уменьшения общности рассуждений можно считать, что $x \in (Ker A)^{\perp}$. Дело в том, что вдоль всякого направления, ортогонального $Ker A$, вспомогательная задача будет имееть обусловленность 1 и не представляет трудоемкости в предположении, что известно $Ker A$ (а, стало быть, и ортогональные направления). Тогда согласно лемме \ref{lemma1} трудоемкость нахождения приемлемого качества решения каждой из задач \eqref{statiiiaux} (линейные слагаемые не влияют на неё) будет определяться числом обусловленности
$$
    \frac{\alpha_i\frac{\lambda_{max} (A^TA)}{\mu_y} + L_x}{\alpha_i\frac{\lambda_{min}^+ (A^TA)}{L_y} + L_x} \leq \frac{L_y \lambda_{max} (A^TA)}{\mu_y \lambda_{min}^{+} (A^TA)},
$$
поскольку для произвольных $a \geq b$ и $c>0$ верно неравенство
$\frac{a+c}{b+c} \leq \frac{a}{b}$. Таким образом, верна оценка \eqref{estimateiii}.

\section{Заключение}\label{section:conclusions}
В данной работе рассмотрен класс гладких седловых задач со структурой. Такие задачи возникают, например, в обработке изображений и при решении различных обратных задач \cite{book:Lan2019}. С помощью оценок ускоренного слайдинга Дж. Лана \cite{book:Lan2019} и результатов о том, как преобразуются свойства гладкости и сильной выпуклости (кривизны графика функции, определяемые экстремальными значениями ее гессиана) при преобразовании Фенхеля–Лежандра \cite{book:Rockafellar}, были получены новые (существенно лучшие) оценки на число вычислений градиента по прямым переменным в нелинейных седловых задачах со структурой. По сути, были получены достаточно общие условия, при которых удается перенести основные результаты для билинейных седловых задач на общий класс нелинейных гладких выпукло-вогнутых седловых задач.

Следует отметить, что по в тексте статьи мы намерено избегали максимальной общности изложения для компактности и удобства восприятия. Тем не менее, отметим, что приведенные в статье результаты можно обобщить на случай более общих оракулов \cite{gasnikov2018book}: стохастических, рандомизированных (в том числе инкрементальных, возникающих при работе с целевыми функционалами вида суммы функций), неполноградиентных, модельных (в том числе с дополнительными композитными членами); вместо евклидовой нормы, можно было рассматривать более общие нормы и прокс-структуры (впрочем, все же не в такой общности, как в не сильно выпуклом случае); наконец можно попробовать рассмотреть более чем два слагаемых в структуре целевого функционала. Также весьма интересен вопрос о том, в какой степени можно перенести полученные результаты на классы негладких седловых задач. Многое упирается в возможность должным образом обобщить лемму \ref{lemma3}. Нам представляется, что на этом пути может быть получено достаточно много новых интересных результатов.

Отметим также, что если в приведенных в статье постановках задач $\dim(y)$ и(или) $\dim(x)$ небольшие, то у использованных в статье ускоренных методов появляются конкуренты в виде методов типа центров тяжести \cite{book:Nemirovski_Ydin_1979,boobk:Bubeck_2015}. В качестве примера можно посмотреть случай $\dim(y)= 2$, разобранный в \cite{paper:Gasnikov_MIPT_2016}.

Несмотря на полученное в работе улучшение известных верхних оценок, мы все же по-прежнему не знаем, оптимальны ли приведенные в этой статье оценки по совокупности критериев: число вычислений градиента по прямым переменным и по двойственным (не только по одному из них)? Недавно появился препринт \cite{paper:lower_bounds_2019}, в котором получены нижние оценки для гладких сильно выпукло-вогнутых седловых задач. Однако при этом была показана достижимость этих оценок только на специальном подклассе задач. В общем случае предложенная в отмеченном препринте \cite{paper:lower_bounds_2019} методика (Section 4.2) приводит к оценкам, аналогичным полученным нами. Однако мы обосновываем указанные верхние оценки для более широкого класса задач, не предполагающих проксимальную дружественность $h$ и $r$.

Авторы выражают благодарность А. С. Немировскому, Ю. Е. Нестерову и Р. Хильдебранду за ценное обсуждение части материала статьи.


\newpage

\section{Приложение}

\subsection{Доказательство леммы \ref{lemma1}}

Ясно, что $g(x)=h^{*}(Ax)$, где $h^{*}$~--- сопряженная функция к $h$. По теореме Демьянова--Данскина $\nabla g(x)=A^{T}y(x)$, где $y(x)=\langle Ax,y(x)\rangle-h(y(x))$
(т.е. $y(x)=\displaystyle\arg\max_{y}\{\langle Ax,y\rangle-h(y)\}$).

Пусть $h(y)$ $\mu_{y}$--сильно выпукла. Тогда в силу выбора $y(x)$ для всяких $x_{1}, x_{2} \in Q_x$:
$$
    \langle Ax_{1},y(x_{2})\rangle-h(y(x_{2}))\leq \langle Ax_{1},y(x_{1})\rangle -h(y(x_1))-\frac{\mu_{y}}{2}\|y(x_{2})-y(x_{1})\|_2^{2},
$$
$$
    \langle Ax_{2},y(x_{1})\rangle-h(y(x_{1}))\leq \langle Ax_{2},y(x_{2})\rangle -h(y(x_{2}))-\frac{\mu_{y}}{2}\|y(x_{1})-y(x_{2})\|_2^{2}.
$$

После сложения этих двух неравенств имеем:
$$
    \langle Ax_{1}-Ax_{2},y(x_{2})-y(x_{1})\rangle\leq -\mu_{y}\|y(x_{2})-y(x_{1})\|_2^{2},
$$
откуда
$$
    \mu_{y}\|y(x_{2})-y(x_{1})\|_2^{2}\leq \langle Ax_{2}-Ax_{1},y(x_{2})-y(x_{1})\rangle\leq \|Ax_{2}-Ax_{1}\|_2\cdot\|y(x_{2})-y(x_{1})\|_2,
$$
т.е. для нормы матрицы $\|A\|_2$
$$
    \|y(x_{2})-y(x_{1})\|_2\leq \frac{\|A\|_2\,\|x_{2}-x_{1}\|_2}{\mu_{y}}.
$$
Поэтому
$$
    \|\nabla g(x_{1})-\nabla g(x_{2})\|_2\leq\|A^{T}\|_2\, \|y(x_{1})-y(x_{2})\|_2\leq\frac{\|A^{T}A\|_2}{\mu_{y}}\|x_{1}-x_{2}\|_2=\frac{\lambda_{max}(A^{T}A)}{\mu_{y}}\|x_{1}-x_{2}\|_2.
$$

Проверим теперь вторую часть утверждения. Пусть теперь $x_{1},x_{2}\in\left(\text{Ker}\, A\right)^{\bot}$.

Хорошо известно, что для сопряженной функции
$$
    h^{*}(x)=\max_{y}\{\langle x,y\rangle-h(y)\}=\langle x,\widehat{y}_{x}\rangle-h(\widehat{y}_{x}),
$$
верно:
$$
    \widehat{y}_{x}\in \partial h^{*}(x)\Longleftrightarrow x\in \partial h(\widehat{y}_{x}),
$$
откуда $(x\rightarrow Ax,\, \widehat{y}_{x}\rightarrow y(x))$
$$
    y(x)\in \partial h^{*}(Ax)\Longleftrightarrow Ax\in \partial h(y(x)).
$$

Тогда имеем:
$$
    \langle\nabla g(x_{1})-\nabla g(x_{2}),x_{1}-x_{2}\rangle=\langle A^{T}y(x_{1})-A^{T}y(x_{2}),x_{1}-x_{2}\rangle=\langle y(x_{1})-y(x_{2}),Ax_{1}-Ax_{2}\rangle=
$$
$$
    =\langle Ax_{1}-Ax_{2},y(x_{1})-y(x_{2})\rangle\geq\{\text{из}\,L_{y} -\text{гладкости}\,h\}\geq\frac{1}{L_{y}}\|Ax_{1}-Ax_{2}\|_2^{2}=
$$
$$
    \frac{1}{L_{y}}\langle A^{T}A(x_{1}-x_{2}),x_{1}-x_{2}\rangle\geq\{\text{из}\,x_{1}-x_{2}\not\in\text{Ker}\, A\,(\text{т.к.}\,x_{1}, x_{2}\in\text{Ker}\,A^{\bot})\}\geq
$$
$$
    \geq\frac{\lambda^{+}_{\min}\left(A^{T}A\right)}{L_{y}}\|x_{1}-x_{2}\|_2^{2},
$$
что обосновывает $\displaystyle\frac{\lambda^{+}_{\min}\left(A^{T}A\right)}{L_{y}}$--сильную выпуклость $g(x)$ при $x\in\left(\text{Ker}\,A\right)^{\bot}$.

\subsection{Доказательство леммы \ref{lemma2}}
Функция $\hat{S}(x,\cdot)$ $\mu_y$-сильно вогнута на $ Q_y $, и $\hat{S}(\cdot,y)$ дифференцируема на $Q_x$. Поэтому по теореме Демьянова--Данскина, для любого $x \in Q_x$
\begin{equation}\label{nabla_g(x)}
    \nabla g(x) = \nabla_x \tilde{S}(x,y^*(x)) = \nabla_x F(x,y^*(x))
\end{equation}

Чтобы доказать, что $g(\cdot)$ имеет $L$--липшицев градиент при $L = L_{xx} + \frac{2 L_{xy}^2}{\mu_y}$, покажем условие Липшица для $y^*(\cdot)$ (функция $y^*$ определена в \eqref{solution:max_problem_g(x)}) с константой $\frac{2L_{xy}}{\mu_y}$.

Ввиду того, что $\hat{S}(x_1, \cdot)$ $\mu_y$-сильно вогнута на $Q_y$, для произвольных $x_1, x_2 \in Q_x$:
\begin{equation}\label{ineq:1}
   \|y^*(x_1) - y^*(x_2)\|_2^2 \leq \frac{2}{\mu_y} \left( \hat{S}(x_1, y^*(x_1)) - \hat{S}(x_1, y^*(x_2)) \right).
\end{equation}

С другой стороны, $\hat{S}(x_2, y^*(x_1)) - \hat{S}(x_2, y^*(x_2)) \leq 0$, т.к. $y^*(x_2)$ доставляет максимальное значение $\hat{S}(x_2,.)$ на $ Q_y $. Имеем
\begin{equation}\label{ineq:2}
\begin{split}
\hat{S}(x_1, y^*(x_1)) - & \hat{S}(x_1, y^*(x_2))  \leq \left( \hat{S}(x_1, y^*(x_1)) - \hat{S}(x_1, y^*(x_2)) \right) - \left( \hat{S}(x_2, y^*(x_1)) - \hat{S}(x_2, y^*(x_2)) \right) \\
&  \hspace{-2em} \stackrel{\text{из} \, \eqref{function:S_hat}}{=} \left( F(x_1, y^*(x_1)) - F(x_1, y^*(x_2)) \right) - \left( F(x_2, y^*(x_1)) - F(x_2, y^*(x_2)) \right) \\
& \hspace{-1em} = \int_0^1 \langle  \nabla_xF(x_1 + t(x_2 - x_1), y^*(x_1)) - \nabla_xF(x_1 + t(x_2 - x_1), y^*(x_2)), x_2 - x_1 \rangle dt \\
&  \hspace{-1em} \leq  \| \nabla_xF(x_1 + t(x_2 - x_1), y^*(x_1)) - \nabla_xF(x_1 + t(x_2 - x_1), y^*(x_1))  \|_2\cdot \|x_2 - x_1\|_2 \\
&   \hspace{-2em} \stackrel{ \text{из}\, \eqref{smooth_F_2}}{\leq} L_{xy} \|y^*(x_1) - y^*(x_2)\|_2\cdot \|x_2 - x_1\|_2.
\end{split}
\end{equation}

Таким образом, из \eqref{ineq:1} и \eqref{ineq:2} вытекает неравенство
\begin{equation}\label{ineq:Lip_y_star}
    \|y^*(x_2) - y^*(x_1)\|_2 \leq \frac{2L_{xy}}{\mu_y}\|x_2 - x_1\|_2
\end{equation}
т.е. функция $y^*(\cdot)$ удовлетворяет условию Липшица с костантой $\frac{2L_{xy}}{\mu_y}$. Далее, из \eqref{nabla_g(x)} получаем
\begin{align*}
\|\nabla g(x_1) -& \nabla g(x_2)\|_2 = \|\nabla_x F(x_1,y^*(x_1))-\nabla_x F(x_2,y^*(x_2))\|_2 \\
&  = \|\nabla_x F(x_1,y^*(x_1)) - \nabla_x F(x_1,y^*(x_2)) + \nabla_x F(x_1,y^*(x_2)) - \nabla_x F(x_2,y^*(x_2))\|_2 \\
&  \leq  \|\nabla_x F(x_1,y^*(x_1)) - \nabla_x F(x_1,y^*(x_2))\|_2 + \|\nabla_x F(x_1,y^*(x_2)) - \nabla_x F(x_2,y^*(x_2))\|_2 \\
&  \hspace{-2em} \stackrel{ \text{из}\, \eqref{smooth_F_1} \, \text{и} \, \eqref{smooth_F_2}}{\leq}  L_{xy}\|y^*(x_1)- y^*(x_2)\|_2 + L_{xx} \|x_2 - x_1\|_2\\
& \hspace{-1em} \stackrel{ \text{из}\, \eqref{ineq:Lip_y_star}}{=} \left( L_{xx} + \frac{2L_{xy}^2}{\mu_y} \right) \|x_2 - x_1\|_2.
\end{align*}

Это означает, что $g(\cdot)$ имеет $L$--липшицев градиент при $L = L_{xx} + \frac{2 L_{xy}^2}{\mu_y}$.

Проверим теперь неравенства из \eqref{eq:lemma2}. Сначала докажем, что для любых $\delta \geq 0$ и $x \in Q_x$ верно
\begin{equation}\label{ineq:inexact_for_g(x)}
    \|\nabla_x \hat{S}(x, \tilde{y}_\delta(x)) - \nabla g(x)\|_2 \leq L_{xy}\sqrt{\frac{2\delta}{\mu_y}}.
\end{equation}

Для всякого $x \in Q_x$ верно $\nabla_x \hat{S}(x, \tilde{y}_\delta(x)) = \nabla_x F(x, \tilde{y}_\delta(x))$. Тогда
\begin{align*}
\|\nabla_x \hat{S}(x, \tilde{y}_\delta(x)) - \nabla g(x)\|_2^2 & = \|\nabla_x F(x, \tilde{y}_\delta(x)) - \nabla_x F(x, y^*(x))\|_2^2 \\
&  \hspace{-1em} \stackrel{ \text{из}\, \eqref{smooth_F_2}}{\leq}  L_{xy}^2 \|y^*(x) - \tilde{y}_\delta(x) \|_2^2 \\
&  \hspace{-1em} \stackrel{ \text{из}\, \eqref{ineq:1}}{\leq}  \frac{2 L_{xy}^2}{\mu_y} \left( \hat{S}(x, y^*(x)) - \hat{S}(x, \tilde{y}_\delta(x))  \right) \\
&  \hspace{-1em} \stackrel{ \text{из}\, \eqref{delta_inexact_solution}}{\leq}  \frac{2\delta L_{xy}^2}{\mu_y},
\end{align*}
что обосновывает неравенство \eqref{ineq:inexact_for_g(x)}.

Теперь в силу $\mu_x$-сильной выпуклости $\hat{S}(\cdot, \tilde{y}_\delta(x))$ на $Q_x$, для произвольных $x, z \in Q_x$ верно
\begin{align*}
    g(z)  \stackrel{ \text{из}\, \eqref{max_problem:g(x)}}{\geq} \hat{S}(z, \tilde{y}_\delta(x))  \geq \hat{S}(x, \tilde{y}_\delta(x))+ \langle \nabla_x\hat{S}(x, \tilde{y}_\delta(x)), z-x \rangle.
\end{align*}

Таким образом,
\begin{equation}
    0 \geq
     \hat{S}(x, \tilde{y}_\delta(x)) - g(z) + \langle \nabla_x\hat{S}(x, \tilde{y}_\delta(x)), z-x \rangle,
\end{equation}
что доказывает левую часть \eqref{eq:lemma2}. Чтобы доказать правую часть \eqref{eq:lemma2}, отметим, что $g$ выпуклая и имеет $L$--липшицев градиент на $Q_x$. Поэтому для произвольных $x, z \in Q_x$ имеем
\begin{align*}
    g(z) & \leq g(x) + \langle \nabla g(x), z-x\rangle + \frac{L}{2}\|z-x\|_2^2\\
    & \hspace{-1em} \stackrel{ \text{из}\, \eqref{delta_inexact_solution}}{\leq}  \hat{S}(x,\tilde{y}_\delta(x)) + \delta + \frac{L}{2} \|z-x\|_2^2 + \langle \nabla g(x), z-x \rangle + \langle \nabla_x\hat{S}(x,\tilde{y}_\delta(x)), x-z \rangle - \\ & - \langle \nabla_x\hat{S}(x,\tilde{y}_\delta(x)), x-z \rangle\\
    & = \hat{S}(x,\tilde{y}_\delta(x)) + \delta + \langle \nabla_x\hat{S}(x,\tilde{y}_\delta(x)), z-x \rangle +\langle \nabla_x\hat{S}(x,\tilde{y}_\delta(x)) - \nabla g(x), x-z \rangle + \frac{L}{2} \|z-x\|_2^2\\
    & \hspace{-1em} \stackrel{ \text{из}\, \eqref{ineq:inexact_for_g(x)}}{\leq} \hat{S}(x,\tilde{y}_\delta(x)) + \delta + \langle \nabla_x\hat{S}(x,\tilde{y}_\delta(x)), z-x \rangle + L_{xy}\sqrt{\frac{2\delta}{\mu_y}}\cdot \|z-x\|_2+ \frac{L}{2}\|z-x\|_2^2.
\end{align*}

Однако
\begin{equation*}
    L_{xy}\sqrt{\frac{2\delta}{\mu_y}}\cdot \|z-x\|_2  \leq  \frac{2\sqrt{\delta}L_{xy}}{\sqrt{\mu_y}} \|z-x\|_2 = 2 \sqrt{\frac{L_{xy}^2}{\mu_y} \|z-x\|_2^2 \cdot \delta} \leq \frac{L_{xy}^2}{\mu_y}\|z-x\|_2^2 + \delta
\end{equation*}
ввиду классического неравенства между средним арифметическим и средним геометрическим. Поэтому
\begin{equation*}
    g(z)  \leq \hat{S}(x,\tilde{y}_\delta(x)) + 2 \delta + \langle \nabla_x\hat{S}(x,\tilde{y}_\delta(x)), z-x \rangle + \frac{L_{xy}^2}{\mu_y}\|z-x\|_2^2 +\frac{L}{2}\|z-x\|_2^2,
\end{equation*}
и поскольку $L = L_{xx} + \frac{2L_{xy}^2}{\mu_y}$, то $\frac{L_{xy}^2}{\mu_y} \leq \frac{L}{2}$ и поэтому
\begin{equation*}
    g(z) \leq \hat{S}(x,\tilde{y}_\delta(x)) + \langle \nabla_x \hat{S}(x,\tilde{y}_\delta(x)), z-x \rangle + 2\delta + L \|z-x\|_2^2.
\end{equation*}

Итак, имеем
\begin{equation}
g(z) - \hat{S}(x, \tilde{y}_\delta(x)) - \langle \nabla_x\hat{S}(x, \tilde{y}_\delta(x)), z-x \rangle \leq L\|z-x\|_2^2 + 2\delta,
\end{equation}
откуда вытекает справедливость левой части неравенства \eqref{eq:lemma2}.

\subsection{Доказательство леммы \ref{lemma3}}
Напомним, что в утверждении леммы \ref{lemma3} рассматривается задача минимизации
	\begin{equation}\label{eq:5}
		\min_{x \in \mathbb{R}^n} P(x) := r(x) + g(x),
	\end{equation}
где функция $r(x)$ $\mu_r$-сильно выпуклая и $L_r$-гладкая для $L_r \geq \mu_r \geq 0$, функция $g(x)$ $\mu_g$-сильно выпуклая и $L_g$-гладкая для $L_g \geq \mu_g \geq 0$, функция $P(x)$ $\mu$-сильно выпукла и $L$-гладкая при $L = L_r + L_g \geq \mu = \mu_r + \mu_g > 0$. Обозначим через $x^*$ искомую точку минмума функционала $P$.

Докажем лемму \ref{lemma3} в предположении, что функция $r(x)$ допускает в произвольной запрошенной точке  $(\delta_r, L_r, \mu_r)$-градиент $\nabla r_{\delta_r}(x)$, а функция $g(x)$ допускает --- $(\delta_g, L_g, \mu_g)$-градиент $\nabla g_{\delta_g}(x)$. Это означает, что для произвольных $x,y\in \mathbb{R}^n$ выполнены неравенства:
		\begin{equation}\label{assump1}
		\frac{\mu_r}{2}\norm{x-y}_2^2 - \delta_r \leq r(x) - r(y) - \langle\nabla r_{\delta_r}(y),x-y \rangle\leq \frac{L_r}{2} \|x-y\|_2^2 + \delta_r,
		\end{equation}
		\begin{equation}\label{assump2}
		\frac{\mu_g}{2}\norm{x-y}_2^2 - \delta_g \leq g(x) - g(y) - \langle\nabla g_{\delta_g}(y),x-y \rangle\leq \frac{L_g}{2} \|x-y\|_2^2 + \delta_g,
		\end{equation}
		где $\delta_r \geq 0$ и $\delta_g \geq 0$.
		
На самом деле для обоснования основных результатов работы достаточно утверждение леммы \ref{lemma3} для менее ограничительной концепции ($\delta, L$)-градиента в предположении сильной выпуклости $g$ и $r$. Поскольку предполагается, что как $r$, так и $g$ допускает неточные значения градиентов в запрашиваемых точках, то для определённости можно положить $L_r \leq L_g$.

Будем применять к рассмотренной задаче \eqref{eq:5} следующий метод, который подразумевает решение вспомогательной подзадачи быстрым градиентным методом при условии неточно заданного ($\delta_g, L_g, \mu_g$)-градиента $g$.
	
	\begin{algorithm}
		\caption{Ускоренный проксимальный градиентный метод c неточными значениями градиентов}
		\begin{algorithmic}[1]
			\STATE {\bf Параметры:} $x^0 \in \mathbb{R}^n$, шаги $\alpha, \beta \in (0,1)$, $\eta > 0$.
			\STATE $y^0 = z^0 = x^0$
			\FOR{$k=0,1,2,\ldots$}
				\STATE $x^k = \alpha z^k + (1-\alpha)y^k$
				\STATE $y^{k+1} \approx \hat{y}^{k+1} := \mathrm{prox}_{\frac{1}{L_r} g(\cdot)}(x^k -  \frac{1}{L_r}\nabla r_{\delta_r}(x^k))$ \label{line:1}\quad($y^{k+1}$~--- приближенное значение данного оператора, найденное посредством решения вспомогательной задачи оптимизации быстрым градиентным методом)
				\STATE $z^{k+1} = \beta z^k + (1-\beta)x^k + \eta (y^{k+1} - x^k)$
			\ENDFOR
		\end{algorithmic}
		\label{alg:1}
	\end{algorithm}

Докажем необходимую вспомогательную оценку для параметров $x^k$ и $y^{k+1}$ при произвольном $x \in \mathbb{R}^n$.

\begin{statement}
		Для всякого $x \in \mathbb{R}^n$ выполнено следующее неравенство:
		\begin{equation}\label{eq:1}
		\begin{split}
		\langle x^k - y^{k+1},x - x^k \rangle\\
		\leq
		\frac{1}{L_r + \mu_g}
		\left[ P(x) - P(y^{k+1}) - \frac{\mu}{4}\norm{x - x^k}_2^2 - \frac{L_r + \mu_g}{4}\norm{y^{k+1} - x^k}_2^2 + 2{\delta_r}\right]\\
		+
		c_1\norm{\hat{y}^{k+1} - y^{k+1}}_2^2,
		\end{split}
		\end{equation}
		где константа $c_1$ определяется следующим выражением:
		\begin{equation*}
			c_1 = 2\left[\frac{L_r}{\mu} + 1\right]
			\left[\frac{L_g^2}{L_r^2} + 1\right].
		\end{equation*}
	\end{statement}

    \begin{proof}  Из определения $\hat{y}^{k+1}$ следует, что
		\begin{equation}\label{Koval1}
		\hat{y}^{k+1} = x^k - \frac{1}{L_r}\nabla r_{\delta_r}(x^k) - \frac{1}{L_r}\nabla g(\hat{y}^{k+1}).
		\end{equation}
		В силу предположения \eqref{assump1} и $\mu_g$-сильной выпуклости функции $g(x)$, имеем
		\begin{align*}
		\langle x^k -& y^{k+1},x - x^k\rangle
		=
		\langle x^k - \hat{y}^{k+1} + \hat{y}^{k+1} - y^{k+1},x - x^k\rangle\\
		&=
		\frac{1}{L_r}\langle \nabla r_{\delta_r}(x^k) + \nabla g(\hat{y}^{k+1}),x - x^k\rangle
		+
		\langle \hat{y}^{k+1} - y^{k+1},x - x^k\rangle\\
		&=
		\frac{1}{L_r}\langle \nabla r_{\delta_r}(x^k) + \nabla g(y^{k+1}),x - x^k\rangle\\
		&+
		\langle \frac{1}{L_r}\left[\nabla g(\hat{y}^{k+1}) - \nabla g(y^{k+1})\right] + \hat{y}^{k+1} - y^{k+1},x - x^k\rangle\\
		&=
		\frac{1}{L_r}\langle \nabla r_{\delta_r}(x^k),x - x^k\rangle
		+
		\frac{1}{L_r}\langle \nabla g(y^{k+1}),x - y^{k+1}\rangle
		+
		\frac{1}{L_r}\langle \nabla g(y^{k+1}),y^{k+1} - x^k\rangle\\
		&+
		\langle \frac{1}{L_r}\left[\nabla g(\hat{y}^{k+1}) - \nabla g(y^{k+1})\right] + \hat{y}^{k+1} - y^{k+1},x - x^k\rangle\\
		&\leq
		\frac{1}{L_r}\left[r(x) - r(x^k) - \frac{\mu_r}{2}\norm{x - x^k}_2^2 + \delta_r\right]
		+
		\frac{1}{L_r}\left[g(x) - g(y^{k+1}) -\frac{\mu_g}{2}\norm{x - y^{k+1}}_2^2 \right]\\
		&+
		\frac{1}{L_r}\langle \nabla g(y^{k+1}), y^{k+1} - x^k\rangle
		+
		\langle \frac{1}{L_r}\left[\nabla g(\hat{y}^{k+1}) - \nabla g(y^{k+1})\right] + \hat{y}^{k+1} - y^{k+1}, x - x^k\rangle.
		\end{align*}
	Далее, применим правую часть неравенства \eqref{assump1} для $r(x)$:
		\begin{equation}
			r(y^{k+1}) \leq r(x^k) + \langle\nabla r_{\delta_r}(x^k),y^{k+1} - x^k\rangle + \frac{L_r}{2}\norm{y^{k+1} - x^k}_2^2 + \delta_r,
		\end{equation}
		откуда следует
		\begin{align*}
			\langle x^k -& y^{k+1},x - x^k\rangle
			\leq
			\frac{1}{L_r}\left[r(x) - r(y^{k+1}) - \frac{\mu_r}{2}\norm{x - x^k}_2^2 + 2{\delta_r}\right]
			+\\
			&+
			\frac{1}{L_r}\left[g(x) - g(y^{k+1}) -\frac{\mu_g}{2}\norm{x - y^{k+1}}_2^2\right]\\
			&+
			\frac{1}{L_r}\langle \nabla r_{\delta_r}(x^k) +  \nabla g(y^{k+1}), y^{k+1} - x^k\rangle
			+
			\frac{1}{2}\norm{y^{k+1} - x^k}_2^2\\
			&+
			\langle \frac{1}{L_r}\left[\nabla g(\hat{y}^{k+1}) - \nabla g(y^{k+1})\right] + \hat{y}^{k+1} - y^{k+1},x - x^k\rangle \\
			&=
			\frac{1}{L_r}\left[ P(x) - P(y^{k+1}) - \frac{\mu_r}{2}\norm{x - x^k}_2^2 - \frac{\mu_g}{2}\norm{x - y^{k+1}}_2^2 - \frac{L_r}{2}\norm{y^{k+1} - x^k}_2^2 + 2{\delta_r} \right]\\
			&+
			\frac{1}{L_r}\langle\nabla r_{\delta_r}(x^k) +  \nabla g(y^{k+1}),y^{k+1} - x^k\rangle
			+
			\norm{y^{k+1} - x^k}_2^2\\
			&+
			\langle \frac{1}{L_r}\left[\nabla g(\hat{y}^{k+1}) - \nabla g(y^{k+1})\right] + \hat{y}^{k+1} - y^{k+1},x - x^k\rangle\\
			&=
			\frac{1}{L_r}\left[ P(x) - P(y^{k+1}) - \frac{\mu_r}{2}\norm{x - x^k}_2^2 - \frac{\mu_g}{2}\norm{x - y^{k+1}}_2^2 - \frac{L_r}{2}\norm{y^{k+1} - x^k}_2^2 + 2{\delta_r}\right]\\
			&+
			\langle \frac{1}{L_r}\left[\nabla g(\hat{y}^{k+1}) - \nabla g(y^{k+1})\right] + \hat{y}^{k+1} - y^{k+1},x^k - y^{k+1}\rangle \\
			&+
			\langle \frac{1}{L_r}\left[\nabla g(\hat{y}^{k+1}) - \nabla g(y^{k+1})\right] + \hat{y}^{k+1} - y^{k+1},x - x^k\rangle \\
		\end{align*}

Применим теперь неравенство Юнга, а также $L_g$-Липшицевость градиента $\nabla g(x)$ $\|\nabla g(y^{k+1}) - \nabla g (\hat{y}^{k+1})\|_2^2 \leq L_g^2 \|\hat{y}^{k+1} - y^{k+1}\|_2^2$:
		\begin{align*}
			\langle x^k -& y^{k+1}, x - x^k\rangle
			\leq
			\frac{1}{L_r}\left[ P(x) - P(y^{k+1}) - \frac{\mu_r}{2}\norm{x - x^k}_2^2 - \frac{\mu_g}{2}\norm{x - y^{k+1}}_2^2 - \frac{L_r}{2}\norm{y^{k+1} - x^k}_2^2\right]\\
			&+ \frac{2\delta_r}{L_r}+\frac{\mu}{4L_r}\norm{x - x^k}_2^2+\frac{L_r + \mu_g}{4L_r}\norm{y^{k+1} - x^k}_2^2\\
			&+
			\left[\frac{L_r}{\mu} + \frac{L_r}{L_r + \mu_g}\right]
			\norm{\frac{1}{L_r}\left[\nabla g(\hat{y}^{k+1}) - \nabla g(y^{k+1})\right] + \hat{y}^{k+1} - y^{k+1}}_2^2\\
			&\leq
			\frac{1}{L_r}\left[ P(x) - P(y^{k+1}) - \frac{\mu_r - \mu_g}{4}\norm{x - x^k}_2^2 - \frac{\mu_g}{2}\norm{x - y^{k+1}}_2^2 - \frac{L_r - \mu_g}{4}\norm{y^{k+1} - x^k}_2^2\right]\\
			&+\frac{2{\delta_r}}{L_r}+2\left[\frac{L_r}{\mu} + \frac{L_r}{L_r + \mu_g}\right]
			\left[\frac{L_g^2}{L_r^2} + 1\right]
			\norm{\hat{y}^{k+1} - y^{k+1}}_2^2.
		\end{align*}
		Наконец, проведём финальные преобразования:
		\begin{align*}
		\langle x^k -& y^{k+1}, x - x^k\rangle =
		\frac{L_r}{L_r + \mu_g} \langle x^k - y^{k+1},x - x^k\rangle + \frac{\mu_g}{L_r + \mu_g}\langle x^k - y^{k+1},x - x^k\rangle\\
		&\leq
		\frac{1}{L_r + \mu_g}
		\left[ P(x) - P(y^{k+1}) - \frac{\mu_r - \mu_g}{4}\norm{x - x^k}_2^2 - \frac{L_r - \mu_g}{4}\norm{y^{k+1} - x^k}_2^2 + 2{\delta_r}\right]\\
		&-
		\frac{\mu_g}{2(L_r + \mu_g)}\norm{x - y^{k+1}}_2^2
		+
		\frac{2L_r}{L_r + \mu_g}\left[\frac{L_r}{\mu} + \frac{L_r}{L_r + \mu_g}\right]
		\left[\frac{L_g^2}{L_r^2} + 1\right]
		\norm{\hat{y}^{k+1} - y^{k+1}}_2^2\\
		&+
		\frac{\mu_g}{2(L_r + \mu_g)}\left[\norm{x - y^{k+1}}_2^2 - \norm{x^k - y^{k+1}}_2^2 - \norm{x - x^k}_2^2\right]\\
		&\leq
		\frac{1}{L_r + \mu_g}
		\left[ P(x) - P(y^{k+1}) - \frac{\mu_r + \mu_g}{4}\norm{x - x^k}_2^2 - \frac{L_r + \mu_g}{4}\norm{y^{k+1} - x^k}_2^2 + 2{\delta_r}\right]\\
		&+
		2\left[\frac{L_r}{\mu} + 1\right]
		\left[\frac{L_g^2}{L_r^2} + 1\right]
		\norm{\hat{y}^{k+1} - y^{k+1}}_2^2.
		\end{align*}
	\end{proof}

\begin{statement}
		Пусть выбраны следующие значения праметров для алгоритма~\ref{alg:1}:
		\begin{equation*}
			\eta = \frac{2(L_r + \mu_g)}{8\alpha(L_r + \mu_g) + (1-\alpha)\mu},
		\end{equation*}
		\begin{equation*}
			\beta = 1 - \frac{\eta\mu}{2(L_r + \mu_g)} = 1 - \frac{\mu}{8\alpha(L_r + \mu_g) + (1-\alpha)\mu},
		\end{equation*}
		\begin{equation*}
			\alpha = \frac{1}{4}\sqrt{\frac{\mu}{L_r + \mu_g}} \leq \frac{1}{4}.
		\end{equation*}
		Тогда выполнено следующее неравенство:
		\begin{equation}\label{eq:6}
		    \begin{split}
		    \norm{z^{k+1} - x^*}_2^2
		+
		c_2
		\left[
		P(y^{k+1}) - P(x^*)
		\right]
		&\leq
		(1-\alpha)
		\left(
		\norm{z^k - x^*}_2^2
		+
		c_2
		\left[
		P(y^k) - P(x^*)
		\right]
		\right)
		\\
		&+
		c_3\left[8c_1\norm{y^{k+1} - \hat{y}^{k+1}}_2^2 - \norm{y^{k+1} - x^k}_2^2  \right]\\
		&+ 4c_2\delta_r.
		\end{split}
		\end{equation}
		где $c_2$ и $c_3$  --- некоторые положительные константы.
\end{statement}


\begin{proof}
Оценим величину $\norm{z^{k+1} - x^*}_2^2$:
		\begin{align*}
			\norm{z^{k+1} - x^*}_2^2
			&=
			\norm{\beta z^k + (1-\beta)x^k - x^* + \eta (y^{k+1} - x^k)}_2^2\\
			&=
			\norm{\beta(z^k - x^*) + (1-\beta)(x^k - x^*)}_2^2
			+
			\eta^2\norm{y^{k+1} - x^k}_2^2
			+\\
			&+
			2\eta\langle \beta z^k + (1-\beta)x^k - x^*,y^{k+1} - x^k\rangle\\
			&\leq
			\beta\norm{z^k - x^*}_2^2
			+
			(1-\beta)\norm{x^k - x^*}_2^2
			+
			\eta^2\norm{y^{k+1} - x^k}_2^2 +\\
			&+
			2\eta\beta\langle z^k - x^k,y^{k+1} - x^k\rangle
			+
			2\eta\langle x^k - x^*,y^{k+1} - x^k\rangle\\
			&\leq
			\beta\norm{z^k - x^*}_2^2
			+
			(1-\beta)\norm{x^k - x^*}_2^2
			+
			\eta^2\norm{y^{k+1} - x^k}_2^2 + \\
			&+
			2\eta\beta\frac{1-\alpha}{\alpha}\ \langle x^k - y^k,y^{k+1} - x^k\rangle
			+
			2\eta\langle x^k - x^*,y^{k+1} - x^k\rangle.
		\end{align*}
		Далее два раза применим неравенство \eqref{eq:1}:
		\begin{align*}
			\norm{z^{k+1} - x^*}_2^2
			&\leq
			\beta\norm{z^k - x^*}_2^2
			+
			(1-\beta)\norm{x^k - x^*}_2^2
			+
			\eta^2\norm{y^{k+1} - x^k}_2^2\\
			&+
			2\beta\frac{\eta}{L_r + \mu_g}\frac{1-\alpha}{\alpha}
			\left[
				P(y^k) - P(y^{k+1}) - \frac{L_r + \mu_g}{4}\norm{y^{k+1} - x^k}_2^2 + 2{\delta_r}
			\right]
			\\
			&+
			2\frac{\eta}{L_r + \mu_g}
			\left[
				P(x^*) - P(y^{k+1}) - \frac{\mu}{4}\norm{x^* - x^k}_2^2 - \frac{L_r + \mu_g}{4}\norm{y^{k+1} - x^k}_2^2 + 2{\delta_r}
			\right]\\
			&+
			2\eta c_1\left[\beta\frac{1-\alpha}{\alpha} + 1\right]\norm{y^{k+1} - \hat{y}^{k+1}}_2^2\\
			&=
			\beta\norm{z^k - x^*}_2^2
			+
			\left[1 - \beta - \frac{\eta\mu}{2(L_r + \mu_g)}\right]\norm{x^k - x^*}_2^2\\
			&+
			\left[
				\eta^2 - \frac{\eta\beta}{4}\frac{1-\alpha}{\alpha} - \frac{\eta}{4}
			\right]\norm{y^{k+1} - x^k}_2^2\\
			&+
			2\beta\frac{\eta}{L_r + \mu_g}\frac{1-\alpha}{\alpha}
			\left[
				P(y^k) - P(y^{k+1})
			\right]
			+
			2\frac{\eta}{L_r + \mu_g}
			\left[
				P(x^*) - P(y^{k+1})
			\right]\\
			&+
			\frac{\eta}{4} \left[\beta\frac{1-\alpha}{\alpha} + 1\right]\left[8c_1\norm{y^{k+1} - \hat{y}^{k+1}}_2^2 - \norm{y^{k+1} - x^k}_2^2  \right]\\
			&+
			2\delta_r\left[
				2\beta\frac{\eta}{L_r + \mu_g}\frac{1-\alpha}{\alpha}
				+
				2\frac{\eta}{L_r + \mu_g}
			\right].
		\end{align*}

		C учётом выбранных значений параметров $\beta$ и $\eta$, а также
		$$
		    c_2 = \frac{2\eta\beta}{\alpha(L_r + \mu_g)}, \quad c_3 = \frac{\eta}{4}\left[\beta\frac{1-\alpha}{\alpha} + 1\right],
		$$
		получаем

		\begin{align*}
		\norm{z^{k+1} - x^*}_2^2
		&\leq
		\beta\norm{z^k - x^*}_2^2\\
		&+
		2\beta\frac{\eta}{L_r + \mu_g}\frac{1-\alpha}{\alpha}
		\left[
		P(y^k) - P(y^{k+1})
		\right]
		+
		2\frac{\eta}{L_r + \mu_g}
		\left[
		P(x^*) - P(y^{k+1})
		\right]+\\
		&+
		c_3\left[8c_1\norm{y^{k+1} - \hat{y}^{k+1}}_2^2 - \norm{y^{k+1} - x^k}_2^2  \right]
		+
		\frac{4\delta_r\eta}{\alpha(L_r + \mu_g)}\\
		&\leq
		\beta\norm{z^k - x^*}_2^2
		+
		\frac{2\beta\eta}{L_r + \mu_g}\frac{1-\alpha}{\alpha}
		\left[
		P(y^k) - P(y^{k+1})
		\right]\\
		&+
		\frac{2\beta\eta}{L_r + \mu_g}
		\left[
		P(x^*) - P(y^{k+1})
		\right]+\\
		&+
		c_3\left[8c_1\norm{y^{k+1} - \hat{y}^{k+1}}_2^2 - \norm{y^{k+1} - x^k}_2^2  \right]
		+
		\frac{4\delta_r\eta}{\alpha(L_r + \mu_g)}\\
		&=
		\beta\norm{z^k - x^*}_2^2
		+
		c_2(1-\alpha)
		\left[
		P(y^k) - P(y^{k+1})
		\right]
		+
		c_2\alpha
		\left[
		P(x^*) - P(y^{k+1})
		\right]+\\
		&+
		c_3\left[8c_1\norm{y^{k+1} - \hat{y}^{k+1}}_2^2 - \norm{y^{k+1} - x^k}_2^2  \right]
		+
		\frac{2c_2\delta_r}{\beta}\\
		&=
		\beta\norm{z^k - x^*}_2^2
		+
		c_2(1-\alpha)
		\left[
		P(y^k) - P(x^*)
		\right]
		+
		c_2
		\left[
		P(x^*) - P(y^{k+1})
		\right]+\\
		&+
		c_3\left[8c_1\norm{y^{k+1} - \hat{y}^{k+1}}_2^2 - \norm{y^{k+1} - x^k}_2^2  \right]
		+
		\frac{2c_2\delta_r}{\beta}.
		\end{align*}
		Используя значение параметра $\alpha$, получаем
		\begin{align*}
			\frac{1}{2} \leq \beta = 1 - \frac{\mu}{2\sqrt{(L_r + \mu_g)\mu} + (1-\alpha)\mu} \leq 1 - \frac{1}{3}\sqrt{\frac{\mu}{L_r + \mu_g}} \leq 1-\alpha,
		\end{align*}
		откуда следует
		\begin{align*}
			\norm{z^{k+1} - x^*}_2^2
			+
			c_2
			\left[
				P(y^{k+1}) - P(x^*)
			\right]
			&\leq
			(1-\alpha)
			\left(
				\norm{z^k - x^*}_2^2
				+
				c_2
				\left[
					P(y^k) - P(x^*)
				\right]
			\right)+
			\\
			&+
			c_3\left[8c_1\norm{y^{k+1} - \hat{y}^{k+1}}_2^2 - \norm{y^{k+1} - x^k}_2^2  \right]
			+
			4c_2\delta_r.
		\end{align*}
	\end{proof}	

Теперь учтем, что вспомогательная задача строки~\ref{line:1} листинга алгоритма~\ref{alg:1} решается быстрым градиентным методом с неточным заданием градиента $g$. Оценим необходимую точность $\delta_g$ градиента $g$ для получения требуемого качества решения задачи по функции.
	
\begin{statement}
		Пусть приближение $y^{k+1}$ прокс-оператора $\hat{y}^{k+1} = \mathrm{prox}_{\frac{1}{L_r} g(\cdot)}(x^k - \frac{1}{L_r}\nabla r_{\delta_r}(x^k))$ (строка~\ref{line:1} листинга алгоритма~\ref{alg:1}) вычисляется быстрым градиентным методом в предположении, что доступен ($\delta_g, \mu_g, L_g$)-градиент $g$ в произвольной запрошенной точке \cite{thesis:devolder2013}. Решаемая задача минимизации при этом имеет вид
		\begin{equation}\label{eq:2}
			\min_{x \in \mathbb{R}^n} g(x) + \frac{L_r}{2}\norm{x^k - \frac{1}{L_r}\nabla r_{\delta_r}(x^k) - x}_2^2,
		\end{equation}
	где $x^k$ --- начальное приближение. Тогда известно \cite{thesis:devolder2013}, что для произвольного $\delta \in (0; 1)$, после
		\begin{equation}\label{eq:3}
			T = O\left(\sqrt{ \frac{L_r + L_g}{L_r + \mu_g}} \log \frac{L_r + L_g}{\delta(L_r + \mu_g)}\right),
		\end{equation}
		итераций указанного метода гарантированно будет выполнено неравенство
		\begin{equation}\label{eq:4}
			\norm{y^{k+1} - \hat{y}^{k+1}}_2^2 \leq \delta\norm{x^k - \hat{y}^{k+1}}_2^2 + c_4\delta_g,
		\end{equation}
		где константа $c_4$ может быть задана выражением
		\begin{equation}
			c_4 = \frac{4\sqrt{L_r + L_g}}{(L_r + \mu_g)\sqrt{L_r + \mu_g}}.
		\end{equation}
\end{statement}	
\begin{proof}
		Отметим, что целевая функция задачи \eqref{eq:2} $(L_r+\mu_g)$-сильно выпуклая и $(L_r + L_g)$-гладкая, а $\hat{y}^{k+1}$ --- точное решение задачи\eqref{eq:2}. Неравенство \eqref{eq:4} следует из соотвествующего результата для быстрого градиентного метода в концепции ($\delta_g, \mu_g, L_g$)-оракула для $g$
		\cite{thesis:devolder2013}.
	\end{proof}

\begin{proof}[Доказательство леммы \ref{lemma3}]
Выбрав в неравенстве \eqref{eq:4} $\delta = \frac{1}{32c_1} \leq \frac{1}{4}$, получим
			\begin{align*}
				\norm{y^{k+1} - \hat{y}^{k+1}}_2^2
				&\leq
				2\delta\left(\norm{x^k - {y}^{k+1}}_2^2 + \norm{y^{k+1} - \hat{y}^{k+1}}_2^2\right)+ c_4\delta_g\\
				&\leq
				2\delta \norm{x^k - {y}^{k+1}}_2^2 + \frac{1}{2}\norm{y^{k+1} - \hat{y}^{k+1}}_2^2 + c_4\delta_g,
			\end{align*}
			откуда следует
			\begin{align*}
			\norm{y^{k+1} - \hat{y}^{k+1}}_2^2
			&\leq
			4\delta\norm{x^k - {y}^{k+1}}_2^2 + 2c_4\delta_g\\
			&\leq
			\frac{1}{8c_1}\norm{x^k - {y}^{k+1}}_2^2 + 2c_4\delta_g.
			\end{align*}
			C учётом доказанных неравенств \eqref{eq:6} означает, что
			\begin{align*}
				\norm{z^{k+1} - x^*}_2^2
				+
				c_2
				\left[
				P(y^{k+1}) - P(x^*)
				\right]
				&\leq
				(1-\alpha)
				\left(
				\norm{z^k - x^*}_2^2
				+
				c_2
				\left[
				P(y^k) - P(x^*)
				\right]
				\right)+\\
				&+
				4c_2\delta_r + 2c_3c_4\delta_g,
			\end{align*}
			откуда после телескопирования имеем
			\begin{align*}
				\norm{z^{k} - x^*}_2^2
				+
				c_2
				\left[
				P(y^{k}) - P(x^*)
				\right]
				\leq
				(1-\alpha)^k &\left(\norm{x^{0} - x^*}_2^2
				+
				c_2
				\left[
				P(x^{0}) - P(x^*)
				\right]\right)+
				\\
				&+
				\frac{4c_2\delta_r + 2c_3c_4\delta_g}{\alpha}.
			\end{align*}
С учётом $\mu$-сильной выпуклости функции $P(x)$ имеем
			\begin{align*}
				P(y^{k}) - P(x^*)
				&\leq
				(1-\alpha)^k\left(1 + \frac{2}{\mu c_2}\right)\left[P(x^0) - P(x^*)\right] + \frac{4\delta_r}{{\alpha}} + \frac{2c_3c_4\delta_g}{c_2\alpha}\\
				&\leq
			2(1-\alpha)^k\left[P(x^0) - P(x^*)\right] +
			\frac{4\delta_r}{{\alpha}} + \frac{2c_3c_4\delta_g}{c_2 \alpha}.
			\end{align*}
			Выбирая число итераций внешнего метода
			\begin{equation}
				k = \frac{1}{\alpha}\log \frac{4(P(x^0) - P(x^*))}{\varepsilon} = O \left(\sqrt{\frac{L_r + \mu_g}{\mu}}\log\frac{1}{\varepsilon}
				\right)
			\end{equation}
			и точность $(\delta_r, L_r, \mu_r)$-градиента $\nabla r_{\delta_r}(x)$
			\begin{equation}
				\delta_r = \frac{\alpha\varepsilon}{16} = O \left(\sqrt{\frac{\mu}{L_r + \mu_g}} \varepsilon\right),
			\end{equation}
			а также $(\delta_g, L_g, \mu_g)$-градиента $\nabla g_{\delta_g}(x)$
			\begin{align}
				\delta_g &= \frac{\alpha c_2\varepsilon}{ 8c_3c_4}
				=
				\frac{\alpha\varepsilon}{8c_4} \frac{2\eta\beta}{\alpha(L_r + \mu_g)} \frac{4\alpha}{\eta[(1-\alpha)\beta + \alpha]}
				\leq
				\frac{\alpha\varepsilon}{c_4(1-\alpha)(L_r + \mu_g)}\nonumber\\
				&=
				\frac{\sqrt{L_r + \mu_g}\alpha\varepsilon}{4(1-\alpha)\sqrt{L_r + L_g}}
				=
				\frac{\sqrt{L_r + \mu_g}\sqrt{\mu}\varepsilon}{16(1-\alpha)\sqrt{L_r + L_g}\sqrt{L_r + \mu_g}}
				\leq
				\frac{\varepsilon}{12}\sqrt{\frac{\mu}{L_r + L_g}}\nonumber\\
				&=
				O\left(\sqrt{\frac{\mu}{L_g + \mu_r}}\varepsilon\right),
			\end{align}
			где в последнем равенстве использовалось предположение $L_r \leq L_g$ и $\alpha \leq \frac{1}{4}$, получаем требуемое качество решения
			\begin{align*}
				P(y^{k}) - P(x^*)
				\leq
				\varepsilon.
			\end{align*}
			При этом количество вызовов $(\delta_r, L_r, \mu_r)$-градиента $\nabla r_{\delta_r}(x)$
			\begin{equation}
				k = O \left(\sqrt{\frac{L_r + \mu_g}{\mu}}\log\frac{1}{\varepsilon}\right),
			\end{equation}
			а количество вызовов $(\delta_g, L_g, \mu_g)$-градиента $\nabla_{\delta_g} g(x)$
			\begin{align*}
				k \times T
				&=
				O \left(\sqrt{\frac{L_r + \mu_g}{\mu}}\log\frac{1}{\varepsilon}\right)\times O \left(\sqrt{ \frac{L_r + L_g}{L_r + \mu_g}} \log \frac{L_r + L_g}{\delta(L_r + \mu_g)}\right)\\
				&=
				\tilde{O} \left(\sqrt{\frac{L_r + L_g}{\mu}}\log \frac{1}{\varepsilon}\right)
				=
				\tilde{O} \left(\sqrt{\frac{L_g + \mu_r}{\mu}}\log \frac{1}{\varepsilon}\right).
			\end{align*}
в силу допущения $L_r \leq L_g$ (данное допущение не существенно ввиду симметричности найденных оценок на $\delta_r$ и $\delta_g$).
	\end{proof}


\bibliography{article}
\bibliographystyle{plain}

\end{document}